%% file: artikel274.tex
\documentclass[10pt]{article}
\usepackage{enumerate, longtable}
\usepackage{epic, graphicx, amsfonts, amsmath, amssymb, theorem,  epsfig}
\usepackage{a4, enumerate, array, longtable}
\usepackage{amsfonts,multirow, amssymb}
\usepackage[utf8]{inputenc}
\usepackage[all]{xy}
\newtheorem{result}{Theorem}
\newtheorem{definition}{Definition}[section]
\newtheorem{theorem}[definition]{Theorem}

\newtheorem{lemma}[definition]{Lemma}
\newtheorem{corollary}[definition]{Corollary}
\newtheorem{proposition}[definition]{Proposition}

\newcommand{\proof}{{\bfseries Proof. } }
\newcommand{\dimv}{\mbox{${\rm \underline{dim}}$}}
\newcommand{\T}[1][\mu]{\mathcal{T}_{#1}}
\newcommand{\pred}{\preccurlyeq}
\newcommand{\arr}[2]{{\renewcommand{\arraystretch}{0.7}\begin{array}{c} #1 \\#2\end{array}}}

\newcommand{\lle}{ \!\left( \!\!\!}
\newcommand{\rri}{\!\!\! \right)\!\!}
\newcommand{\moem}{\!\!\!\!\!}
\newcommand{\arw}{\renewcommand{\arraystretch}{0.6}}

\title{}
\author{}

\begin{document}
\renewcommand{\thefootnote}{\roman{footnote}}
\title{On minimal disjoint degenerations of modules over tame path algebras}

\author{Klaus Bongartz\\Universität Wuppertal\\Germany\thanks{E-mail:bongartz@math.uni-wuppertal.de}\and
\setcounter{footnote}{6} Guido Frank\\Universität Wuppertal\\Germany\thanks{E-mail:frank@math.uni-wuppertal.de} \and Isabel Wolters\\Universität Wuppertal\\Germany}

\maketitle

\begin{abstract}
 For representations of tame quivers the degenerations are controlled by the dimensions of various homomorphism spaces. Furthermore, there is no proper degeneration to an indecomposable. Therefore, up to common direct summands, any minimal degeneration from $M$ to $N$ is induced by a short exact sequence $0 \rightarrow U \rightarrow M \rightarrow V \rightarrow 0$ with indecomposable ends that add up to $N$. We study these 'building blocs' of degenerations and we prove that the codimensions are bounded by two. Therefore, a quiver is Dynkin resp. Euclidean resp. wild iff the codimension of the building blocs is one resp. bounded by two  resp. unbounded. We explain also that for tame quivers the complete classification of all the building blocs is a finite problem that can be solved with the help of a computer.
\end{abstract}

\section*{Introduction}

If an algebraic group acts on a variety, the orbits are partially ordered by inclusion of their closures. Note that there are at least two general methods to determine the orbit closures, namely one method based on Gröbner bases \cite{Gro} and another one proposed recently by Popov \cite{Popov}. But both methods are quite impractical in the special case we are interested in.

This is the action of $G=Gl_{d}$  by conjugation on the variety $Mod_A^d(k)$ of $d$-dimensional representations of an associative finitely generated algebra. The points of this variety are the $A$-module structures on $k^{d}$ and the orbits are in bijection with the isomorphism classes of $d$-dimensional modules. We write $M\leq_{deg}N$ and call $N$ a degeneration of $M$ resp. $M$ a deformation of $N$ iff the orbit to $N$ lies in the closure of the orbit to $M$. Despite a nice representation theoretic characterization obtained by Zwara in  \cite{Extensions}, building on earlier work of Riedtmann in \cite{Riedtmann}, it is in general a hard problem to determine the degeneration order. However, for tame quivers, i.e. quivers whose underlying graph is a Dynkin or an extended Dynkin diagram, the degeneration order on the representations coincides by \cite{Degenerations,Tamedeg} with the partial order $M\leq N$ defined by $[M,X]\leq [N,X]$ for all modules $X$. Here and later on we abbreviate $dim Hom (X,Y)$ by $[X,Y]$ and $dim Ext (X,Y)$ by $[X,Y]^{1}$. Since the indecomposable representations and also the homomorphism spaces between them are known for representations of these quivers ( see \cite{Donovan,Nazarova,zahme,Ringel} ), this gives a  handy description of the degenerations. Moreover, by  \cite{Extensions} the degeneration order is also equivalent to the $\leq_{ext}$- order defined in \cite{Degenerations}. Therefore one has the 'structure theorem' from \cite{Dynkin} for the minimal degenerations  that reduces their classification to the analysis of short exact sequences $0 \rightarrow U \rightarrow M \rightarrow V \rightarrow 0$ with indecomposable end terms such that $M\leq N=U \oplus V$ is a minimal degeneration. In particular, one would like to know the codimension $[N,N]- [M,M]$ of the orbit of $N$ in the orbit closure of $M$ for these 'building blocs'. This is our  main result.
\begin{result}\label{Theo1}
 Let $Q$ be a tame quiver. Then the codimension of any building bloc $M<U\oplus V$ as above is $\leq 2$. If $U$ is preprojective and $V$ is preinjective, the codimension is $1$.
\end{result}

This theorem, proven only by theoretical means, generalizes and simplifies many previous results that were partially based on extensive computer calculations ( \cite{Markolf,Fritzsche,Frank,Wolters} ).

As we will indicate by some figures and tables, there are very many building blocs and their structure remains a mystery to the authors ( see e.g. figure \ref{diaIsabel} in section \ref{defext} or \cite{Wolters} ). In contrast to the cases considered in \cite{Abeasis}, \cite{Riedtmann} or \cite{Degenerations}, it seems to be impossible to get a manageable classification. Nevertheless, as in the preprojective case  treated in \cite{Fritzsche}, there is a certain periodicity that leads to the following somewhat metamathematical statement:

\begin{result}\label{Theo2} The complete classification of all building blocs is a finite problem, that can be solved by computer.
\end{result}

The precise meaning of this statement will be explained in section \ref{classsec}. The first section serves to fix the notation.
Section \ref{reductionsec} introduces our main reduction techniques that allow us to replace one building bloc by another one with modules of smaller dimension. These techniques are contained in the theorems \ref{redtheo} and \ref{divsub}. Subsequently, we turn in the sections \ref{codimsec} and \ref{tubesec} to the study of the codimensions, and we prove theorem 1. In section \ref{concludingsec} we add some remarks and suggestions. Finally, section \ref{tablessec} gives some tables.

This article has its origin in the two dissertations \cite{Frank,Wolters} of the second and the third author. Some of their results were then polished by the first author.

\section{Basic facts and notations}\label{basics}
Throughout this article we work over an algebraically closed field $k$ and consider a tame quiver $Q$.
We assume that the reader is familiar with the basic facts from the representation theory of finite dimensional algebras as given in \cite{Auslander,Ringel} and also with the structure of the module categories over tame quivers as presented in \cite{Donovan,Nazarova,zahme,Ringel}. But we recall shortly the most important facts thereby fixing the notations.

We identify modules over the path algebra $kQ$ and representations of $Q$ over $k$, which will always be of finite dimension. The dimension vector of a $kQ$-module $M$ is denoted by $\dimv(M)$, and for the Auslander-Reiten translations $DTr$ resp. $TrD$ we write shortly $\tau$ resp. $\tau^{-}$. If $M$ and $N$ are modules with the same dimension vectors, one gets from \cite[4.3]{Auslander} the very useful relation $$[M,X]-[N,X]=[\tau^{-}X,M]-[\tau^{-}X,N]$$ for all  modules $X$. 

The Euler form  of $Q$ is denoted by $\langle -,-\rangle$, and its associated quadratic form, the Tits form, by $q$. For modules $X$ and $Y$ one has $$\langle\dimv(X),\dimv(Y)\rangle=[X,Y]-[X,Y]^{1}.$$
Due to Gabriels theorem $kQ$ is representation-finite iff $Q$ is a Dynkin quiver iff $q$ is positive definite. The extended Dynkin ( or Euclidean or affine ) quivers are those quivers where $kQ$ is tame representation-infinite resp. $q$ is positive semi-definite and not positive definite. 

In the Euclidean case, the radical of $q$ is one-dimensional and spanned by the null-root $\delta$ which is a  vector  having positive integral coefficients one of which at least is $1$. 
The defect of a module $X$ is defined by $$\partial(X):=\langle\delta, \dimv(X)\rangle,$$ which is equal to $[E,X]-[E,X]^1$ for any module $E$ with $\dimv(E)=\delta$. 
Another possibility to define the defect uses the coxeter transformation $c$.
 This is the unique endomorphism of $\mathbb{R}^{Q_0}$ that the sends the dimension vector of any non-projective indecomposable to the dimension vector of its $\tau$-translate.
$c$ induces an automorphism of finite order $p(Q)$ on $\mathbb{R}^{Q_0}/\mathbb{R}\delta$ and satisfies
\begin{eqnarray*}
c^{p(Q)}(\dimv(X))=\dimv(X)+\epsilon(Q)\partial(X)\delta
\end{eqnarray*}
for some positive integer $\epsilon(Q)$ depending on $Q$. $p(Q)$ is called the Coxeter number of $Q$, but it should not be confused with the definition of the Coxeter number in Lie Theory.

While for a representation-finite quiver any module is preprojective and preinjective, an indecomposable representation $X$ of an extended Dynkin quiver is either preprojective or regular or preinjective. This is characterized by the defect $\partial (X)$ which is either negative or zero or positive.

Furthermore, a $kQ$-module $X$ has a decomposition $$X=X_P\oplus X_R \oplus X_I,$$ into its preprojective, regular and preinjective parts.

The full subcategory of regular modules breaks up into the direct sum of a $\mathbb{P}^1$-family of so called regular tubes $\T$.
Each regular tube $\T$ is an extension closed abelian subcategory, equivalent to the category $\mathcal{N}(p_\mu)$ of nilpotent representations of the oriented cycle with $p_\mu$ points. The simples $E_1,\; \ldots ,\; E_{p_\mu}$ of this subcategory form a single $\tau$-orbit and their dimension vectors sum up to $\delta$. We always number the simples in the way that $\tau E_{i}=E_{i-1}$ when the indices are read modulo $p_\mu$.
Every indecomposable $R\in\mathcal{T}_\mu$ admits a unique composition series in $\T$. The regular composition factors are then (from the bottom) $S, \tau^{-}S, \ldots, \tau^{-l}S$, for some $l\in \mathbb{N}$. $Soc(R):=S$ is the regular socle, $Top(R):=\tau^{l}S$ the regular top and $l(R):=l+1$ the regular length of $R$. In addition,  the multiplicity of any regular simple $E$ in the regular composition series of
  $R$ is abbreviated by $l_E(R)$ and the module with regular socle $E_{i}$ and regular length $k$ is denoted by $E_{i}(k)$.
The number $p_\mu$ is also called the period of the tube $\mathcal{T}_\mu$. In fact, there are at most three $\mu\in\mathbb{P}^1(k)$ such that $p_\mu\neq 1$. The tubes of period $1$ are called homogeneous. Besides, $p(Q)$ is the least common multiple of the $p_\mu$, $\mu\in\mathbb{P}^1$.

It is important to know that there are no non-zero maps from preinjective to preprojective or regular modules and from regular to preprojective modules.
For indecomposables $R_1$, $R_2$ of the same regular tube $\T$ we have
\[
 [R_1,R_2]=min(l_{Top(R_1)}(R_2),l_{Soc(R_2)}(R_1)).
\]
In particular, if $\T$ is homogeneous, this means $[R_1,R_2]=min(l(R_1),l(R_2)).$

\section{The two reduction techniques}\label{reductionsec}

\subsection{Division by directed summands}
For any module category one defines a preorder $\preceq$ on the set of indecomposables by saying $U \preceq V$ if there is a finite chain of non-zero homomorphisms $f_{i}:V_{i} \rightarrow V_{i+1}$ between indecomposables such that $U=V_{0}$ and $V=V_{n}$. One says that $U$ is a proper predecessor of $V$ if one has $U\preceq V$, but not $V\preceq U$. This preorder is actually a partial order on the preprojective and on the preinjective modules over any path algebra. For extended Dynkin quivers, all modules within the same tube are comparable and modules from different tubes are incomparable.

\begin{definition}
 A decomposition $M=M_{1}\oplus M_{2}\ldots \oplus M_{t}$ of a module $M$ into non-zero direct summands $M_{i}$ is called directed, if for all indecomposable direct summands $U_{i}$ of $M_{i}$ one has that $U_{i}$ is not a  predecessor of $U_{i-1}$ in case $U_{i-1}$ is regular and not a proper predecessor of $U_{i-1}$ in the other cases. The decomposition is disjoint if no indecomposable occurs in two $M_{i}$'s as a direct summand.
\end{definition}

A well-known example of a directed decomposition is the canonical decomposition $M=M_{P}\oplus M_{R}\oplus M_{I}$ of a module into its preprojective, regular and preinjective parts.

\begin{theorem}[Reduction I]\label{redtheo}
Consider a minimal degeneration $M<U\oplus V$ where $U$ is simple projective. Let $M=M_{1}\oplus M_{2}$ be a directed decomposition and let \[ 0\longrightarrow U \stackrel{\epsilon_1\choose \epsilon_2}{\longrightarrow} M=M_1\oplus M_2 
\longrightarrow V \longrightarrow 0,
\] be an exact sequence inducing the minimal degeneration. Then we have:
\begin{enumerate}
 \item $C= Coker (\epsilon_{1})$ is indecomposable.
\item $\Delta:=Codim(U\oplus V,M)-Codim(U\oplus C,M_1)$ is equal to 
\[ [V,V]-[C,C]-[V,M_2]+[V,M_2]^1-[M_2,M_2]^1-[M_2,M_1].\]
\item If $M_{1}$ and $M_{2}$ are disjoint, then the induced degeneration $M_{1}<U \oplus C$ is minimal again.
\end{enumerate}

\end{theorem}
\proof The proof of (a) is essentially same as the proof of theorem 1 in \cite{Fritzsche}. We only sketch the main steps. Since $U$ is simple projective, 
there is the following commutative diagram:
\begin{eqnarray}\label{reddiagr}
\parbox{5cm}{
\xymatrix{
            &                      &0\ar[d]                       &0\ar[d]           &\\
            &                      &M_2\ar@{=}[r]\ar[d]           &M_2       \ar[d]   &\\
0\ar[r]     &U\ar[r]\ar@{=}[d]     &M_1\oplus M_2\ar[r]\ar[d]     &V\ar[r]\ar[d]     &0\\
0\ar[r]     &U\ar[r]               &M_1\ar[r]\ar[d]               &C\ar[r]\ar[d]     &0\\
            &                      &0                             &0                 &
}}\end{eqnarray}
If we assume that $C$ is not indecomposable, we can write $C=C_1\oplus C_2$ with $C_1$ indecomposable and $C_2\neq 0$. This induces two commutative diagrams, the second one by applying the snake lemma:

\[
\xymatrix{
            &                      &0\ar[d]                       &0\ar[d]           &\\
 0\ar[r]    &U\ar[r]\ar@{=}[d]     &M_1'\ar[r]\ar[d]              &C_1\ar[r]\ar[d]   &0\\
0\ar[r]     &U\ar[r]               &M_1\ar[r]\ar[d]               &C\ar[r]\ar[d]     &0\\
            &                      &C_2\ar@{=}[r]\ar[d]           &C_2    \ar[d]     &\\
            &                      &0                             &0                 &\\
&&&&
}\qquad
\xymatrix{
            &                      &                              &0\ar[d]           &\\
            &0\ar[d]               &                              &C_2\ar[d]   &\\
0\ar[r]     &M_2\ar[r]\ar[d]       &V\ar[r]\ar@{=}[d]             &C\ar[r]\ar[d]     &0\\
0\ar[r]     &M_2'\ar[r]\ar[d]      &V\ar[r]                       &C_1\ar[r]\ar[d]     &0\\
            &C_2 \ar[d]            &                              &0                 &\\
&0&&&
}
\]
In particular, there are exact sequences
\[\begin{array}{cl}
 (i) &0 \rightarrow U \rightarrow M_1'\rightarrow C_1\rightarrow 0\\
 (ii) &0 \rightarrow M_1' \rightarrow M_1\rightarrow C_2\rightarrow 0\\
(iii) &0 \rightarrow M_2' \rightarrow V\rightarrow C_1\rightarrow 0 \mbox{ and }\\
(iv) &0 \rightarrow M_2 \rightarrow M_2'\rightarrow C_2\rightarrow 0.
\end{array}\]
Using these sequences, 
 it follows that $M_1\oplus M_2 \leq M_1'\oplus M_2' \leq U\oplus V$. 
Because $C_1\neq 0$, the minimality of $M<U\oplus V$ forces $M'\cong M$. On the other hand, the sequence $(iv)$ 
  implies the existence of an indecomposable direct summand $X$ of $M_2'$ that also occurs in $M_1$. Furthermore, if $V$ is regular, $X$ is as a direct summands of $M_1$ non-regular. It satisfies $X\prec V\preceq X$, which is absurd.  Accordingly, $C$ must be  indecomposable.

(b) Due to the projectivity of $U$ and the exact sequences $0\rightarrow U\rightarrow M\rightarrow V\rightarrow 0$ resp. $0\rightarrow M_2\rightarrow V\rightarrow C\rightarrow 0$, we obtain
\begin{eqnarray*}
  [M,M_2]&=&[V,M_2]+[U,M_2]-[V,M_2]^1+[M,M_2]^1 \mbox { resp. }\\
\; 0&=&[U,M_2]-[U,V]+[U,C].
\end{eqnarray*}
Since the decomposition $M=M_1\oplus M_2$ is directed, it follows that $[M_1,M_2]^1=0$. 
From there, we get
\begin{eqnarray*}
\Delta&=&[U\oplus V, U\oplus V]-[M,M]-[U\oplus C,U\oplus C]+[M_1,M_1]\\
&=&[V,V]-[C,C] + [U,V]-[U,C] -[M,M_2]-[M_2,M_1]\\
&=&[V,V]-[C,C]-[V,M_2]+[V,M_2]^1-[M_2,M_2]^1-[M_2,M_1].
\end{eqnarray*}

(c) Assume that $M_1<U\oplus C$ is not minimal. Then there exists a module  $N$  with $M_1<N<U\oplus C$ such that $0\rightarrow U\rightarrow N\rightarrow C\rightarrow 0$ is minimal. $N$ can be decomposed as follows:
\begin{itemize}
\item[\textbullet] $N_1$ contains all indecomposable direct summands $Y$ of $N$ such that there exists an indecomposable direct summand $X$ of $M_1$ with $Y\preceq X$.
\item[\textbullet] $N_2$ consists of the remaining direct summands.
\end{itemize}
$0<[M_1,M_1]\leq[N,M_1]=[N_1,M_1]$ guarantees that $N_1$ is non-zero. Moreover, there is an injection $N_2\hookrightarrow C$, which induces the following commutative pullback diagram
\begin{eqnarray*}
\xymatrix{
0\ar[r]   &M_2\ar[r]\ar@{=}[d]   &P\ar[r]\ar@{^{(}->}[d]  &N_2\ar[r]\ar@{^{(}->}[d]    &0\\
0\ar[r]   &M_2\ar[r]             &V\ar[r]                     &C\ar[r]&0}
\end{eqnarray*}
In particular, $P$ is non-zero and degenerates to $M_2\oplus N_2$.

We claim that $M< N_1\oplus P< U\oplus V$, which is a contradiction to the minimality of $M<U\oplus V$. Let $T$ be an indecomposable module. If $T$ is a predecessor of some direct summand of $M_1$, then we have $0=[M_2,T]=[N_2,T]$. Thus, also $[P,T]$ vanishes. We obtain
\[
 [M,T]=[M_1,T]\leq [N_1\oplus N_2,T]=[N_1\oplus P,T]\leq[U\oplus C,T]\leq[U\oplus V,T].
 \]
If $T$ is not injective and there is no such summand in $M_1$,  then $0=[\tau^{-}T,M_1]=[\tau^{-}T,N_1]$. The injections $M_2\hookrightarrow P \hookrightarrow V$ imply
\begin{eqnarray*}
[U\oplus V,T]-[N_1\oplus P,T]&=& [\tau^{-}T,U\oplus V]-[\tau^{-}T,N_1\oplus P]\\
&\geq& [\tau^{-}T,V]-[\tau^{-}T,P]\;\geq \;0 \quad\mbox{ and }\\
\;[N_1\oplus P,T]-[M,T]&=& [\tau^{-}T,N_1\oplus P]-[\tau^{-}T,M]\\
&=& [\tau^{-}T,P]-[\tau^{-}T,M_2]\;\geq\; 0.
\end{eqnarray*}
Finally, for injective $T$ the equality of the dimension vectors leads to
\[[U\oplus V,T]=[N_1\oplus P,T]=[M,T].\] 
 
If $M$ and $N_1\oplus P$ were isomorphic, then the construction of the decomposition $N=N_1\oplus N_2$ would imply that $P$ is the direct sum of $M_2$ and certain direct summands of $M_1$. But  we have $M_2\hookrightarrow P$, whence $P=M_2$. This leads to $N_2=0$ and $M_1=N$, a contradiction. On the other hand, assumed that $N_1\oplus P\cong U\oplus V$, there are the following possibilities:\\
$N_1\cong U$: Then we would have $N_2\cong C$, which violates $N<U\oplus C$.\\
$N_1\cong V$: This would force $P\cong U$, whence $M_2\cong U$. This is impossible, since $M<U\oplus V$.\\
So $M_1<U\oplus C$ must be minimal. $\hfill \square$\\

A simple consequence of theorem \ref{redtheo} is the following

\begin{corollary}\label{preprojdef}
If $U$ is simple projective and $V$ is preinjective and $M<U\oplus V$ is minimal, then $\partial(M_P)>\partial(U)$. Dually, also  $\partial(M_I)<\partial(V)$ holds. 
\end{corollary}
\proof This is clear for $M_P=0$. Otherwise, choose in theorem \ref{redtheo} $M_1=M_P$ and $M_2= M_R \oplus M_I$. Now consider the sequence
$
 0 \rightarrow U \rightarrow M_P \rightarrow C \rightarrow 0.
$
We have $V\preccurlyeq C$, whence $C$ is again preinjective. This shows $\partial(M_P)=\partial(U\oplus C)>\partial(U).\hfill \square$

\subsection{Division by certain submodules}
We will use some definitions and constructions from \cite[2.1]{Dynkin} that we recall briefly. Given two modules $E$ and $M$, the module  $Q$ is called the generic quotient of $M$ by $E$ if there is a mono $f:E \longrightarrow M$ with cokernel $Q$ and if all cokernels of monos from $E$ to $M$ are degenerations of $Q$. Clearly, a generic quotient need not exist. But it does if there are only finitely many isomorphism classes of cokernels of monos, because the constructible set of all cokernels is irreducible.

Similarly, $M$ is the generic extension of $Q$ by $E$ if $M$ is an extension having all other extensions as a degeneration. Again, a generic extension need not exist, but it does provided there are only finitely many isomorphism classes of extensions.

 Of course, the two operations are in general not inverse to each other.

\begin{theorem}[Reduction II]\label{divsub}
 let $A$ be a finite dimensional algebra with modules $E,M,M',Q,Q'$. Introduce the following conditions:
\begin{enumerate}[i)]
 \item $[E,M]=[E,M']$ and $[Q,E]^{1}-[Q,E]=[Q',E]^{1}-[Q',E]$.
\item $Q$ is the generic quotient of $M$ by $E$ and $M$ is the generic extension of $Q$ by $E$.
\item $Q'$ is the generic quotient of $M'$ by $E$ and $M'$ is the generic extension of $Q'$ by $E$.
\item For all $L$ with $Q\leq L \leq Q'$ there is a generic extension with $E$ and $Q$ is the only quotient of $M$ by $E$ that degenerates to $Q'$.
\item All $N$ with $M\leq N \leq M'$ have a generic quotient by $E$ and $M$ is the only extension of $Q$ by $E$ that degenerates to $M'$.
\end{enumerate}
Then we have:
\begin{enumerate}[a)]
\item If i), ii) and iii) hold, then  $M\leq M'$ is equivalent to $Q\leq Q'$ and we have $Codim (M,M')= Codim (Q,Q')$.
\item If i), ii), iii) and iv) hold,  then $Q<Q'$ is minimal provided $M<M'$ is so. 
\item Finally, if i), ii), iii) and v) hold, then $M<M'$ is minimal provided $Q<Q'$ is minimal.
\end{enumerate}
\end{theorem}
\proof First suppose that only i), ii) and iii) hold. By \cite[theorem 1]{Dynkin}  $M\leq M'$ is equivalent to $Q\leq Q'$.

In the proof of \cite[theorem 1]{Dynkin} one defines a variety $c^{-1}\mathcal{Q}$ and  two smooth morphisms $\lambda: c^{-1}\mathcal{Q} \longrightarrow \overline{O(q)}$ and $\rho: c^{-1}\mathcal{Q} \longrightarrow \overline{O(m)}$ of relative dimensions $l$ and $r$ to the orbit closures of $M$ and $Q$ ( see \cite{Hartshorne} ).  By construction and by the condition iii) we have $ \lambda^{-1}(O(q') \subseteq \rho^{-1}( \overline{O(m')})$ and $ \rho^{-1}(O(m'))\subseteq \lambda^{-1}(\overline{O(q')})$.  Similarly one obtains $ \lambda^{-1}(O(q) \subseteq \rho^{-1}( \overline{O(m)})$ and $ \rho^{-1}(O(m))\subseteq \lambda^{-1}(\overline{O(q)})$ from condition ii). Using the well-known formulas for the dimensions of the fibres of a flat morphism ( see \cite{Hartshorne} again )  we obtain
$dim \,O(m) + r= dim \,O(q)+l$ and $dim \,O(m') + r= dim \,O(q')+l.$  The wanted relation $Codim (M,M')= Codim (Q,Q')$ follows immediately.

 Now suppose that iv) holds in addition. Let $M<M'$ be minimal and take a module $Q''$ with $Q<Q''<Q'$. Then the generic extension $M''$ of $Q''$ by $E$ satisfies $M\leq M'' \leq M'$ as one  sees by looking at an appropriate bundle of cocycles ( \cite[section 2.1]{Dynkin} ). By minimality, we have $M''=M$ or $M''=M'$. In the first case, $M$ has two non-isomorphic quotients by $E$ that degenerate to $Q'$, and in the second case $Q'$ is not the generic quotient of $M'$ by $E$. So $Q<Q'$ is also minimal.

Finally assume that i), ii), iii) and v) hold. Let $Q<Q'$ be minimal and take a module $M''$ with $M<M''<M'$. Then the generic quotient $Q''$ of $M''$ by $E$ satisfies $Q\leq Q'' \leq Q'$  because of \cite[theorem 2.4]{Degenerations}. But $Q''=Q$ contradicts the assumption that $M$ is the only extension of $Q$ by $E$ degenerating to $M'$, and $Q''=Q'$ implies that $M'$ is not the generic extension of $Q'$ by $E$.$\hfill \square$\\

\section{The codimension of the building blocs}\label{codimsec}
\subsection{Deformations between regular modules}\label{regularssec}

Here, we assume that $U$ and $V$ are regular in the same tube $\mathcal{T}$ with regular simples $S_{1},S_{2}, \ldots ,S_{p}$. Then $\mathcal{T}$ is equivalent to the nilpotent representations of an oriented cycle with $p$ points, and this case was already treated by Kempken in \cite{Kempken} a long time ago and in a different language. Therefore we give here an independent proof thereby also specifying which deformations are extensions. In the following proof $S[m]$ denotes the  indecomposable of regular length $m$ with regular top $S$, in particular $S[0]=0$.
\begin{lemma}\label{codimreg}
Let $U=S_i[k]$ and $V=S_j[l]$ be indecomposables in $\mathcal{T}$ and let $r$ be the minimal length of an indecomposable module $W$ with $Top(W)=Top(V)$ and $Soc(W)=\tau^{-}Top(U)$. Then the partially ordered sets 
	$$\mathcal{S}(V,U)= \left\{m | m \in \mathbb{N}\,,\,l \geq r+mp > l-k  \right\}$$ and
	$$\mathcal{E}(V,U)= \left\{M | M \not\simeq U\oplus V ,\mbox{ there is an exact sequence }0 \to U \to M \to V \to 0  \right\}$$ 
	are in bijection under the order-reversing map $m \mapsto S_i[k+r+mp]\oplus S_j[l-r-mp]$.\\
	For the unique minimal element $M$ in $\mathcal{E}(V,U)$ we have
	\begin{equation*}
	\text{Codim}(U \oplus V, M) = \begin{cases}
								~~	2 &\text{ , for } l(U) \geq l(V) \text{ and } Top(U) = Top(V) \\
								~~	2 &\text{ , for } l(U) < l(V) \text{ and } Soc(U) = Soc(V) \\
								~~  1 &\text{ , otherwise.}
						 	 				 					\end{cases}
\end{equation*} 
\label{reg}
\end{lemma} 

\proof
First we show that $S_i[k+r+mp]\oplus S_j[l-r-mp]$ belongs to $\mathcal{E}(V,U)$ for $m \in \mathcal{S}(V,U)$. Set $M_1=S_i[k+r+mp]$. Then we have $Soc(M_1)=U$ and $l(M_1)>l(U)$. So there is a proper monomorphism $\varepsilon_1 : U \to M_1$. Similarly we have $Top(M_1)= Top(V)$, $l(M_1)>l(V)$ and a proper epimorphism $\pi_1 : M_1 \to V$. For $l = r+mp$ there is an obvious exact sequence $0 \to U \to M_1 \to V \to 0$. In the other case $K = \ker \pi_1$ is a proper submodule of $U$ with canonical projection $\varepsilon_2 : U \to U / K $. Set $M_2=U / K$ and $\varepsilon = \left(
												\begin{smallmatrix}\varepsilon_1 \\
 																					\varepsilon_2
 											 	\end{smallmatrix}\right)$
and look at the exact sequence $$0 \to U \xrightarrow{\varepsilon} M_1 \oplus M_2 \xrightarrow{\pi=\left( \begin{smallmatrix}\pi_1 &	\pi_2 \end{smallmatrix}\right)} C \to 0.$$ By construction, $Top(\varepsilon_2)$ is an isomorphism, whence also $Top(\pi_1)$. Counting lengths we see that $C=S_j[l]$ and $M_2=S_j[-r-mp]$.\\
The injectivity of the map is obvious. To see that the map is surjective we take a non-split exact sequence $0 \to U \xrightarrow{\varepsilon} M \xrightarrow{\pi} V \to 0$. The induced exact sequence $0 \to Soc(U) \to Soc(M) \to Soc(V)$ shows $l(Soc(M)) \leq l(Soc(U))+ l(Soc(V)) = 2$. So $M$ is indecomposable or the direct sum $M_1 \oplus M_2$ of two indecomposables. In the first case we have $M=S_i[k+l]=S_i[k+r+mp]$ with $r+mp=l$. So $M$ is in the image of the map. For $M=M_1 \oplus M_2$ we assume $l(M_1) \geq l(M_2)$ and decompose $\varepsilon$ and $\pi$. Now, the induced sequences on the socles and on the tops are again exact by a length argument. Since the sequence does not split,  one of the $\varepsilon_i$ is a proper mono and one a proper epi. The same holds for  the $\pi_k$. Because of $l(M_1) \geq l(M_2)$ we conclude that $\varepsilon_1$ is mono and $\pi_1$ is epi.
So we have $M_1=S_i[k+r+mp]$ and $M_2=S_j[-r-mp]$ for some $m$ with $l > r+mp > l-k$.

Now we take $m < m+1$ in $\mathcal{S}$ and show that $X<Y$ holds for
$$X=S_i[k+r+(m+1)p] \oplus S_j[l-r-(m+1)p]$$ and $$Y=S_i[k+r+mp] \oplus S_j[l-r-mp].$$ Take any indecomposable $Z$ with $l(Z) \leq k+r+mp$. Then the image of any $f ~:~ Z \to S_i[k+r+(m+1)p]$ has length $\leq k+r+mp$. So it factors through $S_i[k+r+mp] \hookrightarrow S_i[k+r+(m+1)p]$ and we have $[Z,S_i[k+r+mp]] = [Z,S_i[k+r+(m+1)p]]$. Because of $[Z,S_j[l-r-mp]] \geq [Z,S_j[l-r-(m+1)p]]$ the inequality $[Z,X] \leq [Z,Y]$ holds. If $l(Z) > k+r+mp$, we have $[Z,Y] = l_{Top(Z)}(Y)= l_{Top(Z)}(X) \geq [Z,X]$.

Finally, to derive the codimension formula we can assume that $k \geq l$ up to duality. The minimal element in $\mathcal{E}(V,U)$ is then given by $S_i[k+r] \oplus S_j[k-r] = M_1 \oplus M_2$. We calculate the codimension 
\begin{align*}
c=&~[U \oplus V,U \oplus V] - [M,M]\\
=&\left( [U \oplus V,U \oplus V] - [U \oplus V,M] \right) + \left( [U \oplus V,M] - [M,M] \right).
\end{align*}
Since $M_1$ has maximal length and $0 \to U \to M \to V \to 0$ is exact we have\\
$[U \oplus V,M_1]=[M,M_1]=l_{Soc(M_1)}(M)$. The surjection $M_1 \xrightarrow{\pi_1} V$ induces an isomorphism $Hom(V,M_2) \simeq Hom(M_1,M_2)$ since any $f: M_1 \to M_2$ has $\ker \pi_1$ in its kernel because of $l(\ker f) > l(\ker \pi_1)$. The surjection $U \xrightarrow{\varepsilon_2} M_2$ also induces an isomorphism $Hom(U,M_2) \simeq Hom(M_2,M_2)$. So we get $[U \oplus V,M] - [M,M]=0$.\\
The inclusion $U \hookrightarrow M_1$ gives  $[U,U]=[U,M_1]$ and  $[V,U]=[V,M_1]$. We always have $[V,V] - [V,M_2]=l_{Top(V)}(V)-l_{Top(V)}(M_2) \leq 1$ and $[V,V] - [V,M_2]=1$ because the identity does not factor through the inclusion $M_2 \hookrightarrow V$. So we get
\begin{align*}
c=&~1+ [U ,V] - [U,M_2]\\
=&~1 + l_{Top(U)}(V) - l_{Top(U)}(M_2) = 1 + l_{Top(U)}(V / M_2) \\
=&~1 + l_{Top(U)}(M_1 / U).
\end{align*}
The wanted formula is now obvious.$\hfill\square$

\subsection{An inductive codimension formula}\label{indcodimsec}

By lemma \ref{codimreg} we already know the codimension of the  building blocs where $U$ and $V$ are regular.
So, up to duality and tilting,  we can assume from now on that $U=P(x)$ is the simple projective corresponding to  the only sink $x$ in $Q$.

If $V$ is non-regular, theorem \ref{redtheo} leads to a very useful inductive formula  to analyze the codimension of a minimal degeneration $M<U\oplus V$.  Suppose that $M$ has a disjoint directed decomposition $M=M_1\oplus M_2$. We consider the commutative diagram 
\begin{eqnarray*}
\xymatrix{
            &                      &0\ar[d]                       &0\ar[d]           &\\
            &                      &M_2\ar@{=}[r]\ar[d]           &M_2       \ar[d]   &\\
0\ar[r]     &U\ar[r]\ar@{=}[d]     &M_1\oplus M_2\ar[r]\ar[d]     &V\ar[r]\ar[d]     &0\\
0\ar[r]     &U\ar[r]               &M_1\ar[r]\ar[d]               &C\ar[r]\ar[d]     &0\\
            &                      &0                             &0                 &
}\end{eqnarray*}
By theorem \ref{redtheo} $M_1<U\oplus C$ is minimal again. Since $V$ is non-regular, we have $[V,V]=1$ and $[V,M_2]=0$.
If $V$ is preinjective, then also $C$ is, which implies $[C,C]=1$. Otherwise, $V$ is preprojective and $C$ can belong to any connected component of the Auslander-Reiten quiver. For preprojective or preinjective $C$ we  obtain immediately $[C,C]=1$. In the case where $C$ is regular, i.e. $M_1$ is preprojective of the same defect as $U$, lemma 3 of \cite{Fritzsche} insures  that $\dimv(C)\leq \delta$, i.e. $[C,C]=1$.
Therefore we have
\[
\Delta:=Codim(U\oplus V,M)-Codim(U\oplus C,M_1)=[V,M_2]^1-[M_2,M_2]^1\geq 0.
\]

Apart from that, dualization delivers the minimal de\-ge\-ne\-ra\-tion $DM=DM_2\oplus DM_1 < DV\oplus DU$ of $kQ^{op}$-modules.
 Choosing a slice $\mathcal{S}$ in the Auslander-Reiten quiver of $kQ^{op}$ with $DV$ as source we can define a tilting module $T:=\bigoplus_{X\in\mathcal{S}}X$ whose endomorphism algebra is a path algebra with the same underlying graph as $Q$. We set $F:=Hom(T,\_)$.

Notice, if $Q$ is Euclidean, the defect behaves in the following way under application of the functor $F\circ D$. If $X$ is some $kQ$-module such that $DX\in\mathcal{T}(T)$, there is
\begin{eqnarray*}
\partial(FDX)=\left\{
 \begin{array}{ll}
-\partial(X),&\mbox{ if } V\in\mathcal{I}\\
\phantom{-}\partial(X),&\mbox{ if } V\in \mathcal{P}.
\end{array}
\right.
\end{eqnarray*}
Furthermore, non-homogenenous tubes of period $p_\mu$ are mapped into non-homoge\-neous tubes of period $p_\mu$. 

Tilting the above situation via $F$ yields 
a new minimal de\-ge\-ne\-ra\-tion $FDM<FDV\oplus FDU$.  $FDV$ is simple projective, $FDU$ is non-regular, $FDM=FDM_2\oplus FDM_1$ is a directed decomposition and there is the following commutative diagram:
\begin{eqnarray*}
\xymatrix{
            &                      &0\ar[d]                       &0\ar[d]           &\\
            &                      &FDM_1\ar@{=}[r]\ar[d]           &FDM_1       \ar[d]   &\\
0\ar[r]     &FDV\ar[r]\ar@{=}[d]     &FDM_2\oplus FDM_1\ar[r]\ar[d]     &FDU\ar[r]\ar[d]     &0\\
0\ar[r]     &FDV\ar[r]               &FDM_2\ar[r]\ar[d]               &L\ar[r]\ar[d]     &0\\
            &                      &0                             &0                 &
}\end{eqnarray*}
By theorem \ref{redtheo} (c) $FDM_2<FDV\oplus L$ is minimal of codimension
\begin{eqnarray*}
&& \hspace{-1.25cm}Codim(FDV\oplus L,FDM_2)=[FDV\oplus L,FDV\oplus L]-[FDM_2,FDM_2]\\
&=&1+[FDV\oplus L,L]-[FDM_2,L]+[FDM_2,L]-[FDM_2,FDM_2]\\
&=&1+[FDV\oplus L,L]-[FDM_2,L]+[FDM_2,FDV]^1-[FDM_2,FDM_2]^1\\
&=&1+[L,L]^1+[V,M_2]^1-[M_2,M_2]^1.
\end{eqnarray*}
This implies $\Delta=Codim(FDV\oplus L,FDM_2)-1-[L,L]^1$, so we have proved the following theorem.

\begin{theorem}{\bf (Inductive codimension formula)} 
\label{decofcodim}
 Under the above assumptions $Codim(U\oplus V,M)$ is equal to
\[ Codim(U\oplus C,M_1)+Codim(FDV\oplus L,FDM_2)-1-[L,L]^1.\]
\end{theorem}

 The power of this result is best illustrated if $V$ is preinjective. Then $L$ is also preinjective, whence $[L,L]^1=0$. 
Thus, inductive application of theorem \ref{decofcodim}  delivers a decomposition of the codimension in the following sense: We can write
\[
 Codim(U\oplus V , M) = 1+\sum_{i=1}^r (Codim(U_i\oplus V_i, M_i)-1)
\]
where  $M_i<U_i\oplus V_i$ is a minimal de\-ge\-ne\-ra\-tion of $kQ_i$-modules such that $U_i$ is projective simple, $V_i$ is preinjective indecomposable, $M_i$ has no proper disjoint directed decomposition and the underlying unoriented graphs of the quivers $Q_{i}$ all coincide.  Thus to show that the codimension is $1$ in general one only has to analyze the following special cases:
\begin{eqnarray*}
\begin{array}{cl}
 \mbox{ $(i)$ }& M_i=X^t \mbox{ where }X \mbox{ is preprojective indecomposable, }\\
 \mbox{ $(ii)$ }& M_i=X^t \mbox{ where }X \mbox{ is preinjective indecomposable or }\\
 \mbox{ $(iii)$ }& M_i=M_\mu\in \T \mbox{ for some }\mu\in\mathbb{P}^1.
\end{array}
\end{eqnarray*}
 
The first and second case are dual to each other. They will be treated in the next subsection by a direct method. The third case is much more complicated. It will be dealt with in the next section.

\subsection{A special case of codimension 1}

Here we prove in a special case directly that the codimension is one.

\begin{lemma}\label{codimzugross}
 Suppose $M<U\oplus V$ is a minimal degeneration. 
\begin{enumerate}
 \item If $Codim(U\oplus V,M) > 1+[V,V]^{1}$, there is an indecomposable direct summand $X$ of $M$ with $[X,M]<[X,U\oplus V]$.
\item If $[X,M]<[X,U\oplus V]$ or $[M,X]<[U\oplus V,X]$ for some indecomposable direct summand $X$ of $M$, we have  $M\oplus X=M'\oplus Z <U'\oplus V'\oplus Z < U\oplus V\oplus X$ for a minimal degeneration $M'<U'\oplus V'$ and some module $Z$. Here $X$ occurs with multiplicity $\geq 2$ in $M'$.
\end{enumerate}
\end{lemma}
\proof a) The minimal degeneration comes from an exact sequence $0\longrightarrow U \longrightarrow M \longrightarrow V \longrightarrow 0$. The induced exact sequence $$0 \rightarrow Hom(V,V) \rightarrow Hom(M,V) \rightarrow Hom(U,V) \rightarrow Ext(V,V)$$ implies $[U\oplus V,U\oplus V]-[M,U\oplus V] \leq 1+ [V,V]^{1}$. So if the codimension is strictly greater than $1+[V,V]^{1}$, we get from 
\begin{eqnarray*}
 Codim(U\oplus V, M)&=&[U\oplus V,U\oplus V]-[M,U\oplus V] + [M,U\oplus V]-[M,M]\\
&\leq&1+[V,V]^1+[M,U\oplus V]-[M,M] 
\end{eqnarray*}
that there is a direct summand $X$ with $[X,M]<[X,U\oplus V]$. 

b) Up to duality we can assume $[X,M]<[X,U\oplus V]$  for some indecomposable direct summand $X$ of $M$. By theorem 4 in \cite{Dynkin} $M\oplus X < U\oplus V \oplus X$ is no longer minimal and we may insert a minimal degeneration $L$ of $M\oplus X$ in between. Again by theorem 4 in \cite{Dynkin} we have $M\oplus X=M'\oplus Z < L = U'\oplus V'\oplus Z$. Here $X$ cannot be a direct summand of $Z$ because the original degeneration is minimal. Thus $X$ occurs with multiplicity $\geq 2$ in $M'$.$\hfill \square$\\

\begin{proposition}\label{indanfang}
Let $M$ be a preprojective or regular semisimple module such that $M<U\oplus V$ is  minimal and $U$ simple projective. Suppose  $M=\bigoplus_{i=1}^sM_i^{r_i}$ is the decomposition into indecomposables. If $End(V)=k$ and $[M_i,M_j]=0$ for $i\neq j$, then $Codim(U\oplus V ,M)\leq 1+[V,V]^1$.
\end{proposition}
\proof The proof is essentially the same as that of proposition 5 in \cite{Dynkin}. Only the beginning has to be modified.
Assume $Codim(U\oplus V,M) > [V,V]^1+1$. By the lemma above there is  some index $i$ that satisfies $[M_i,U\oplus V]>[M_i,M]=r_i$. Without loss of generality $i=1$ may be assumed. This allows us to choose a set of linearly independent  homomorphisms 
\begin{itemize}
\item[\textbullet] $f_{1,1},f_{1,2},\ldots ,f_{1,r_1+1}$ in $Hom(M_1,V)$ resp.
\item[\textbullet]  $f_{k,1},f_{k,2},\ldots ,f_{k,r_k}$ in $Hom(M_k,V)$, $2\leq k\leq s$.
\end{itemize} 
We take two homomorphisms $g_1:M_1^{r_1}\rightarrow V, \quad m\mapsto (f_{1,1},\ldots, f_{1,r_1})(m)$, $g_2:M_1^{n_1}\rightarrow V, \quad m\mapsto (f_{1,2},\ldots, f_{1,r_1+1})(m)$.
For $(a,b)\in k^2\setminus \{(0,0)\}$ we define $f_{(a,b)}
:M=\bigoplus_{k=1}^sM_k^{r_k}\rightarrow V$  as follows:
\[
 f_{(a,b)}(m)=\left\{
\begin{array}{ll}
 (ag_1+bg_2)(m),& ,m\in M_1^{r_1}\\
(f_{k,1},\ldots ,f_{k,r_k})(m), &,m\in M_k^{r_k}, k\geq 2.
\end{array}\right.
\]
 Here it is  convenient to denote the $l$-th copy of $M_k$ in $M$ by $M_{k,l}$ and components of a map $h$ starting or ending there by $h_{k,l}$.
We consider the exact sequence 
\[0\rightarrow K\xrightarrow{(g_{ k,l})} M\xrightarrow{f_{(a,b)}} V\]
 where $K$ is the kernel of $f_{(a,b)}$. $K$ is preprojective. This is trivial for preprojective $M$. If $M$ is regular semi-simple and $K$ would not be  preprojective, then $K$ would contain one of the regular simple summands $M$ contradicting  the choice of the $f_{i,j}$. The remainig part of the proof works as in \cite{Dynkin}.$\hfill \square$\\

\subsection{Deformations with multiple preprojective summands}

The results obtained in this subsection are of independent interest, but they are also essential in the proof of proposition \ref{codimVreg}.
\begin{proposition} \label{multpreprojsum}
 Let $M < U \oplus V$ be a minimal degeneration with a directed decomposition $M=M_{1}\oplus M_{2}$ such that $M_{1}$ and $M_{2}$ are  preprojective and disjoint.
Suppose $M_{1}=M_{1}'\oplus X^{n}$ with $n\geq 2$ and $X$ an indecomposable $\preceq$-maximal direct summand of $M_{1}$ not occurring any more in $M_{1}'$. Then we have:
\begin{enumerate}
 \item $V$ is preprojective, $Q$ an extended Dynkin quiver and $\dimv(M_{1}\oplus X)= \delta + \dimv( U)$.
\item $n=2$ and all indecomposable direct summands of $M$ not isomorphic to $X$ have multiplicity one.
\item $\partial X=-1$.
\end{enumerate}

\end{proposition}
\proof a) and b): For any $i=1,\ldots n$ look at the exact sequence $0 \rightarrow U \rightarrow M_{1}'\oplus X^{i} \rightarrow C_{i} \rightarrow 0$ induced by the directed decomposition $M=M_{1}'\oplus X^{i}\oplus X^{n-i}\oplus M_{2}$. By part a) of theorem \ref{redtheo}, all the dimension vectors $\dimv(C_i)$ are roots of the Tits form $q$. But also $\dimv(C_0)=\dimv(M_{1}')-\dimv(U)$ is a root. This follows as above for $M_{1}'\neq 0$ and it is trivial for $M_{1}'=0$. 

Of course, we have $\dimv(C_i)=\dimv(C_0) + i \dimv( X)$ and 
$$q(\dimv(C_i))=q(\dimv(C_0)) + 2i \langle\dimv(C_0),\dimv(X)\rangle + i^{2},$$
 because $\dimv(X)$ is a real root of $q$.
If $q(\dimv(C_0))=0$, then $\dimv(C_0)=m \delta
$ and $\dimv(C_2)$ is not a root. If $q(\dimv(C_0))=1$, a short calculation with the above equations for $i=1,2$ delivers $q(\dimv(C_2))=2q(\dimv(C_1)) +1$, whence $\dimv(C_1)$ belongs to the radical of $q$.  Thus $Q$ is an extended Dynkin-quiver and $\dimv(C_{3})$ is not a root, which implies $n=2$. Furthermore, $V$ is preprojective because its homomorphic image $C_{2}$ is so having the same defect as $X$. So we get from lemma 3 in \cite{Fritzsche} that $\dimv(C_1)=\delta$. Finally, take any decomposition $M=\oplus_{i=1}^{m}U_{i}$ into indecomposables such that $U_{k}\preceq U_{l}$ implies $k\leq l$. Then the dimension vectors of the quotients $\oplus_{i=1}^{j}U_{i}/U$ are strictly increasing with $j$. Thus the value $\delta$ can be hitted only once. This implies that the other indecomposable direct summands of $M$ occur with multiplicity one only.

c) By part c) of theorem \ref{redtheo} the exact sequence $0 \rightarrow U \rightarrow M_{1}'\oplus X^{2} \rightarrow C_{2} \rightarrow 0$ induces again a minimal degeneration between preprojectives. This time, $X$ is $\preceq$-maximal in $M$. Using duality and tilting, we can even reduce to the case $X^{2}<U \oplus V$ with preprojective indecomposables $U,V$. By proposition \ref{indanfang}, the codimension is one. We obtain
$1=[U\oplus V, U\oplus V]-[U\oplus V,X^{2}]+[U\oplus V,X^{2}]-[X^{2},X^{2}]= 1 + [U\oplus V,X^{2}]-[X^{2},X^{2}],$
whence $[U,X]=2$. But for any homogeneous simple $H$ we have $\partial X= [H,X]-[H,X]^{1}=[X,X]-[X,X]^{1}-[U,X]+[U,X]^{1}=-1$. $\hfill \square$\\

As an example with multiplicities  we choose an $\tilde{E}_{8}$-quiver in the orientation where $U=P(3)$ is the only simple projective and take $M=\oplus_{i=0}^{10}TrD^{i}P(8) \oplus TrD^{5}P(8)$ and $V=TrD^{15}P(3)$.

\subsection{V preinjective and M preprojective}

\begin{proposition}\label{codimDynkin} Let $M$ be a preprojective module having the direct sum of a preprojective indecomposable $U$and a preinjective indecomposable $V$ as a minimal degeneration. Then the codimension is one. This holds in particular for all minimal disjoint degenerations for representations of Dynkin-quivers.
 
\end{proposition}
\proof It is clear that in the decomposition of the codimension, described after theorem \ref{decofcodim}, the third case does not occur. Thus we are done by proposition \ref{indanfang}.$\hfill \square$

\subsection{V regular}

In this subsection $V$ is a regular indecomposable. From $M=M_{P}\oplus M_{R}\leq U\oplus V$ one gets easily $M_{R}\subseteq V$ by using $[X,M_{R}]\leq [X,V]$ for the regular socle $X$ of $M_{R}$ and for $X=M_{R}$.

\begin{proposition}\label{codimVreg} Let $M_{P}\oplus R \leq U\oplus V$ be a  degeneration. Set $C=V/R$.  Then the following holds:
\begin{enumerate}
\item $M_{P}\oplus R \leq U\oplus V$ is minimal iff $M_{P}\leq U\oplus C$ is minimal. Furthermore we have $Codim(U\oplus V,M_{P}\oplus R)= Codim (U\oplus C,M_{P})$.
\item  $\dimv(C)
\leq \delta$.
\item $Codim(U\oplus C,M_{P})=1+[C,C]^{1}$. Thus the codimension is bounded by two.
\end{enumerate}
\end{proposition}
\proof a) We check that all assumptions of theorem \ref{divsub} are satisfied for $E=R$, $M=M_{P}\oplus R$, $M'= U\oplus V$, $Q=M_P$ and $Q'=U\oplus C$. We have $[R,M]=[R,R]=[R,V]=[R,M']$. $[Q,E]^{1}-[Q,E]=[Q',E]^{1}-[Q',E]$ holds because our algebra is hereditary and the dimension vectors of $Q$ and $Q'$ coincide. Also $Q$ is the only quotient of $M$ by $E$ and $M$ the only extension of $Q$ by $E$. Furthermore, $Q'$ is the only quotient and $M'$ the generic extension. Finally, there are generic quotients and generic extensions since there are always only finitely many candidates. Part a) is now a consequence of theorem \ref{divsub}.

b) Suppose $\dimv(C) > \delta$. Then there are two short exact sequences $0 \longrightarrow R_{1} \longrightarrow C
\longrightarrow R_{2} \longrightarrow 0$ and $0 \longrightarrow R_{2}' \longrightarrow C
\longrightarrow R_{1}' \longrightarrow 0$ with $\dimv(R_{1})=\dimv(R_{1}')=\delta$ and $R_{2}\neq 0 \neq R_{2}'.$ Look at the pullback $0 \longrightarrow U \longrightarrow M' \longrightarrow R_{1} \longrightarrow 0$ of the given sequence $0 \rightarrow U \rightarrow M=M_{P} \rightarrow C \rightarrow 0$ by the inclusion $R_{1}\longrightarrow C$. We claim that $M < M'\oplus R_{2}' <U\oplus C$ holds. 

 For an indecomposable preprojective $T$ we have $[M,T]\leq [M',T]+[R_{2},T] =[M',T]\leq [U,T]=[U,T]+[R_{1},T].$ For regular $T$ we have $0=[T,M]\leq [T,M']+[T,R_{2}']=[T,R_{2}']\leq [T,C]=[T,U]+[T,C].$ Finally, for preinjective $T$ we get $[T,M]=[T,M'\oplus R_{2}']=[T,U\oplus C]=0$. Since all three modules have the same dimension vectors, the Auslander-Reiten-formula mentioned in section \ref{basics} implies $M\leq M'\oplus R_{2}' \leq U\oplus C$. All inequalities are strict and this is a contradiction to the minimality.

c) Applying $Hom(-,C)$ to $0 \rightarrow U \rightarrow M_{P} \rightarrow C \rightarrow 0$ one gets the inequality $$Codim (U\oplus C,M_{P}) \geq [U\oplus C,U\oplus C]- [M_{P},C]= 1+[C,C]^{1}.$$ Suppose the inequality is strict. By lemma \ref{codimzugross} we have $M_{P}\oplus X= M'\oplus Z
< U'\oplus V' \oplus Z < U\oplus C \oplus X$ where $X^{2}$ is a direct summand of $M'$. Let $M'=M_{1}'\oplus X^{2}\oplus M_{2}'$ be a directed decomposition. Then we know from proposition \ref{multpreprojsum} that this decomposition is disjoint and that  $\dimv(M_{1}'\oplus X)
=\delta + \dimv(U')$. From $\dimv(M_{1}'\oplus X \oplus M_{2}'\oplus Z) = \dimv(U\oplus C)$ we get $\dimv(C)=\delta$, $U=U'$  and  $M_{2}'=Z=0$. Thus $X$ is a $\preceq$-maximal direct summand of $M$ with defect $-1$ by proposition \ref{multpreprojsum}. In the induced sequence $0 \longrightarrow U \longrightarrow M_{1}' \longrightarrow C' \longrightarrow 0$ the right end $C'$ is preinjective and $Codim( U\oplus C',M_{1}')=1$ by proposition \ref{codimDynkin}. Now from part b) of theorem \ref{redtheo} we obtain for the difference $\Delta$ of the codimensions
$$\Delta= [C,C]-[C',C']-[C,X]+[C,X]^{1}-[X,X]^{1}-[X,M_1']=[C,X]^{1}=-\partial X =1. $$ This is a contradiction.$\hfill \square$\\

\subsection{V preprojective}\label{codimpre}
\begin{proposition} If $M<U\oplus V$ is a minimal degeneration between preprojectives, the codimension is $\leq 2$.
 
\end{proposition}
 \proof We make an induction on the number of indecomposable direct summands of $M$. If $M$ is a power of an indecomposable, the codimension is $1$ by proposition \ref{indanfang}. For the induction we choose a disjoint directed decomposition $M=M_{1}\oplus M_{2}$ with $M_{2}=X^{n}$ for some indecomposable $X$. We look at the induced minimal degenerations $M_1<U\oplus C$ and $FDM_{2}<FDV\oplus L$ as in section \ref{indcodimsec}. Here we always have $Codim(FDV\oplus L,FDM_{2})
-1 - [L,L]^{1}=0$. For preprojective $L$ this is true by proposition \ref{indanfang}, for regular $L$ by proposition \ref{codimVreg} c) and finally for preinjective $L$  by proposition \ref{codimDynkin}. Thus we obtain $Codim(U\oplus V,M)=Codim(U\oplus C,M_{1})$ from theorem \ref{decofcodim}. If $C$ is preprojective our claim follows by induction. For regular $C$ we can use proposition \ref{codimVreg} c) once more and for preinjective $C$ proposition \ref{codimDynkin}.$\hfill\square$

\section{The building blocs with $\mathbf{M}$ in one tube }\label{tubesec}

Throughout this section, $U=P(x)$ is the only simple projective and $V$ is preinjective.
To complete the proof of theorem \ref{Theo1} it remains to consider minimal degenerations $M<U\oplus V$ such that  $M=M_\mu$ comes from a single tube $\T$.

\subsection{A finite test for degenerations}

But first, we derive a finite test criterion for degenerations $M\leq U\oplus V$ that holds for general $M$. Let $M$ be a module with the same dimension vector as $U\oplus V$. We write $M=M_P\oplus \bigoplus_{\mu\in\mathbb{P}^1}M_\mu\oplus M_I$  where $M_P$ is the preprojective part, $M_\mu\in\T$ and $M_I$ is the preinjective part of $M$.

In the following proposition, $d(Q)$ denotes the diameter of $Q$, i.e. the number of edges in the longest unoriented path in $Q$ without cycles. Furthermore, for two indecomposables $X$, $Y$  belonging to the same connected component of the Auslander-Reiten quiver, the length of a shortest path leading from $X$ to $Y$ is abbreviated by $d(X,Y)$. If  $X\not \pred Y$, we set $d(X,Y):=-\infty$.

\begin{proposition}{\bf(Degeneration Test)}\label{enttest}
Under the above assumptions the mo\-dule $M$ degenerates to $U\oplus V$ if and only if it satisfies the following conditions:
\begin{enumerate}[(i)]
 \item $[U,T]-[M_P,T]\geq 0$ for any indecomposable preprojective $T$ such that there exists some direct summand $X$ of $M_P$ with $d(X,T)\leq 2(p(Q)+d(Q))$.
\item $[T,V]-[T,M_I]\geq 0$ for any indecomposable preinjective $T$ such that there exists some direct summand $X$ of $M_I$ with $d(T,X)\leq 2(p(Q)+d(Q))$.
\item  $0<d:=\partial(V)-\partial(M_{I})$ and $[U,E]\geq[M,E]$ for each regular simple $E$. 
\end{enumerate}In that case, each $M_{\mu}$ has at most $d$ indecomposable direct summands.
\end{proposition}
\proof We begin with the necessity of these conditions. 
Since $M$ degenerates to $U\oplus V$, it follows  for any preprojective resp. preinjective $T$ that
\begin{eqnarray*}
 [U,T]-[M_P,T]&=&[U\oplus V, T]-[M,T]\geq 0 \quad\mbox{ resp.}\\\;
[T,V]-[T,M_I]&=&[T,U\oplus V]-[T,M]=[U\oplus V, \tau T]-[M,\tau T]\geq 0.
\end{eqnarray*}
 Hence, the conditions (i) and (ii) hold. For (iii), let $M_\mu^1, \ldots ,  M_\mu^t$
denote the indecomposable direct summands of $M_\mu$. Suppose $E\in\T$ is regular simple, then we get
\[
 [M_\mu^i,E]=\left\{
\begin{array}{ll}
 1& Top(M_\mu^i)=E\\
0& Top(M_\mu^i)\neq E.
\end{array}\right.
\]
This implies (iii).

Reversely, we have to verify the inequality $[U\oplus V,T]-[M,T]\geq 0$ for all indecomposable non-injective $T$. If $T$ is preprojective but no successor of any indecomposable direct summand  of $M_P$, the assertion is clear. Supposed  $T$ is preprojective such that there is no indecomposable direct summand $X$ of $M_P$ with $d(X,T)\leq 2(p(Q)+d(Q))$, we choose $k$ minimal with $d(X,T)\leq 2(k+1)(p(Q)+d(Q))$ for at least one of these $X$. Then  $\tau^{kp(Q)}T$ satisfies $d(X,\tau^{kp(Q)}T)\leq 2(p(Q)+d(Q))$ and by the minimality of $k$ we obtain
\begin{eqnarray*}
[M,T]=[M_P,T]&=&\langle\dimv(M_P),\dimv(T)\rangle\\
&=&\langle\dimv(M_P),\dimv(\tau^{kp(Q)}T)-k\epsilon(Q)\partial(T)\delta\rangle\\
&=&[M_P,\tau^{kp(Q)}T]+k\epsilon(Q)\partial(T)\partial(M_P)\\
&\leq& [U,\tau^{kp(Q)}T]+k\epsilon(Q)\partial(T)\partial(U)=[U,T]=[U\oplus V,T].
\end{eqnarray*}
The dual argument works for preinjective indecomposables. The subsequent application of the Auslander-Reiten formula yields the desired inequality for preinjective but not injective $T$.\\
Finally, let $T$ be a regular module. By our assumptions we have $[U,E]-[M,E]\geq 0$ for all simple regular modules. Let $T$ be an extension of $T'$ and a regular simple $E$. Then we have $[U,T]=[U,T']+[U,E]$ and $[M,T]\leq [M,T']+[M,E]$. This implies $[U\oplus V,T]=[U,T]\geq [M,T]$ by induction on the regular length.$\hfill \square$\\

If  $M$ comes from a single tube $\mathcal{T}$, the criterion simplifies to

\begin{corollary}{\bf (Degeneration test for regular modules)}\label{regenttest}
 For a module $M\in\mathcal{T}$ the following conditions are equivalent:
\begin{enumerate}[(i)]
\item $M<U\oplus V$.
\item $\dimv(M)=\dimv(U\oplus V)$ and  $[U,E]\geq[M,E]$ for each regular simple $E\in \mathcal{T}$.
\end{enumerate}
\end{corollary}

\subsection{Basic facts about tubes}
Throughout the remaining part of this section we consider a minimal degeneration $M=M_{1}\oplus M_{2}\ldots \oplus M_{s}<U \oplus V$ such that  all indecomposable direct summands $M_{i}$ of $M$ belong to the same tube $\mathcal{T}$ with period $p$ and regular simples $E_{1},E_{2},\ldots E_{p}$. Recall, our convention is $\tau E_{i}=E_{i-1}$. We have $\partial V=[\oplus_{i=1}^{p}E_{i},U\oplus V]=\sum_{i=1}^{p}[E_{i},V]$ and $\sum_{i=1}^{p}[E_{i},M]=s$, because $[E_{i},M]$ is the number of indecomposable summands $M_{j}$ with regular socle $E_{i}$.  So $M<U\oplus V$ implies $s\leq \partial V$ with equality holding iff $[E_{i},V]=[E_{i},M]$ for all $i$. Dually we also have $s \leq -\partial U$ with equality iff $[U,E_{i}]=[M,E_{i}]$ for all $i$. Here $[M,E_{i}]$ is the number of modules $M_{j}$ with regular top $E_{i}$.

   We need the following lemma about the 'biology' of the regular simples.
\begin{lemma}\label{gleichverteilung}
 Let $\mathcal{T}$ be  a tube of period $p\geq 2$. Then we have:
\begin{enumerate}
 \item For any preprojective indecomposable $U'$ one has $|[U',E]-[U',F]| \leq 1$ for all regular simples $E$,$F$ in $\mathcal{T}$. In particular, one has $[U',E]\leq 3$ because the defects are bounded by $6$. The dual statement is also true.
\item Let $E_{1},E_{2}, \ldots E_{p}$ be the regular simples. Suppose $\partial(V) \geq p+1$. Remove one of the simples $E_{i}$ from the list. Then there is a point where all the remaining simples vanish or there is an arrow not ending in $x$ that is represented by an isomorphism in all the remaining simples.
\end{enumerate}
\end{lemma}
\proof Part a) is trivial for type $\tilde{A}_{n}$ and $\tilde{D}_{4}$. For $\tilde{D}_{n}$ with $n\geq 5$ one verifies that an indecomposable $E$ with $dim E(y)=2$ and $dimE(z)=0$ at  the two branching points either has not defect $0$ or it has a proper submodule of defect $0$. For the remaining cases one has to consult the table \ref{regsimplestable}. This remarkable equi-distribution can be explained to some extent by the wings in \cite{Ringel} and by the distribution of the roots in the Auslander-Reiten quiver of a Dynkin-quiver, but we take it here as a 'biological' fact.

To prove in part b) that an arrow $\alpha$ is represented by an isomorphism in a regular simple $E$, one only has to check  $dim E(y)= dim E(z)$ for the two end points $y,z$ of the arrow. Namely, $Ext(E,E)=0$ implies that the orbit of $E$ is open, so that by the irreducibility of the module variety it meets the non-empty set where the arrow is represented by an isomorphism. Thus part b) can also be verified by looking at the table \ref{regsimplestable}.$\hfill \square$\\
 
The next result is an easy consequence of this lemma. 

\begin{lemma}\label{regkurze} Let $M<U\oplus V$ be a minimal degeneration. Then the following is  true:
\begin{enumerate}
\item If we have $[U,E]>[M,E]$ for some regular simple $E$, then $E$ occurs only in the regular top of $M$ and any indecomposable summand of $M$ has length $\leq p$. 
\item If we have $[U,E]>[M,E]=0$ and $\partial V \geq p+1$, then $E$ does not occur as a regular composition factor and the codimension is one. 
\end{enumerate}
 The dual statements are  also true.
\end{lemma}
\proof Suppose $[U,E]>[M,E]$, but there is an indecomposable  direct summand $X$ of $M$ having a proper submodule $X'$ with regular top $E$. Then we have by our degeneration test \ref{regenttest} $M=M'\oplus X < M'\oplus X'\oplus X/X' <U \oplus V$ contradicting the minimality of $M<U\oplus V$.

 For $[M,E]=0$ it follows that $E$ does not occur as a composition factor.  If all the other simples vanish at the same point, then $M$ and $U\oplus V$ are obviously representations of a Dynkin-quiver. The claim holds by proposition \ref{codimDynkin}. If one arrow not ending in $x$ is represented by an isomorphism in all the remaining simples, the same is true for $M$ and $U\oplus V$. So we can shrink this arrow and end up again with representations of a Dynkin quiver.$\hfill \square$

\begin{lemma}\label{gleichlang2}
Let $M<U\oplus V$ be minimal. If $E_k(l)$ and $E_q(r)$  are direct summands of $M$, then we have $q=k+\alpha\; mod \; p$ and $r=l+\beta-\alpha\; mod \; p$  for some $0\leq \alpha, \beta\leq p$. In particular, two indecomposable summands satisfy $|l(M_i)-l(M_j)|\leq p$.
\end{lemma}
\proof
Suppose for instance, there occur two indecomposable direct summands $E_k(l)$ and $E_{k+\alpha}(l+p+t-\alpha)$ such that  $t>0$, $ 0\leq \alpha \leq p$. Then we can replace $E_k(l)$ by $E_k(l+t)$  and $E_{k+\alpha}(l+p+t-\alpha)$ by $E_{k+\alpha}(l+p-\alpha)$ to obtain a degeneration between $M$ and $U\oplus V$. For the other possibilities we  proceed similarly.$\hfill \square$\\

The statement of the above lemma can be visualized as follows. If $E_k(l)$ is a direct summand of $M$, then the remaining summands of $M$ do not lie outside the square of the picture down below (resp. the part of the square that actually exists in  case of $l< p_\mu$).
\begin{center}
\setlength{\unitlength}{0.7cm}
\begin{picture}(8,6)
\put(3,3){\line(1,1){2}}    \put(3,3){\line(1,-1){2}}
\put(5,1){\line(1,1){2}}   \put(7,3){\line(-1,1){2}}
\dashline{0.2}[30](4,0)(5,1)
\dashline{0.2}[30](6,0)(5,1)
\dashline{0.2}[30](0,0)(3,3)
\dashline{0.2}[30](0,6)(3,3)
\dashline{0.2}[30](7,3)(8,4)
\dashline{0.2}[30](7,3)(8,2)
\dashline{0.2}[30](5,5)(4,6)
\dashline{0.2}[30](5,5)(6,6)
 \put(3,3){\circle*{0.2}}     \put(7,3){\circle*{0.2}}
\put(5,5){\circle*{0.2}}      \put(5,1){\circle*{0.2}}
\put(3.4,2.8){$\scriptstyle E_k(l)$}
\put(7.4,2.8){$\scriptstyle E_k(l)$}
\put(5.4,0.8){$\scriptstyle E_k(l-p)$}
\put(5.4,4.8){$\scriptstyle E_k(l+p)$}
\end{picture}
\end{center}

The last result in this subsection says what $Codim(U\oplus V,M)=1$ means for the indecomposable direct summands of $M$.

\begin{lemma}\label{gleichlang}
  Let $\mathcal{T}$ be a tube of rank $p$ with simples $E_{1},E_{2},\ldots E_{p}$. Let $M<U\oplus V$ be a minimal degeneration. Decompose $M=\oplus_{i=1}^{p}\oplus _{j=1}^{n_{i}}M_{ij}$ into indecomposables such that the regular top of $M_{ij}$ is $E_{i}$. Then we have:
\begin{enumerate}
 \item If $[U,E_{i}]=n_{i}$ for all $i$ and all $M_{ij}$ have the same regular length, then the codimension is $1$.
\item Reversely, if $\partial V \geq p$, if all $E_{i}$ occur in $M$ as composition factors and if the codimension is $1$, then we have $[U,E_{i}]=n_{i}$ for all $i$ and all $M_{ij}$ have the same regular length. 
\item In particular, if $\mathcal{T}$ is homogeneous and the codimension is one, we obtain $M=R^{\partial(V)}$ for some indecomposable $R$.
\end{enumerate}

 \end{lemma}
\proof Of course, we have 
$$[U,E_{i}]\geq [M,E_{i}]=n_{i} \mbox{ and } [U,M_{rs}]\geq [M,M_{rs}]$$ for all $i,r,s$. For the codimension we get
\begin{eqnarray*}
 Codim(U\oplus V,M)&=&[U\oplus V,U\oplus V]-[U\oplus V,M]+[U\oplus V,M]-[M,M]\\
     &=&1+ \sum_{i,j}([U,M_{ij}]-[M,M_{ij}]).
\end{eqnarray*}
Thus the codimension is $1$ iff $[U,M_{rs}]=[M,M_{rs}]$ holds for all $M_{rs}$. But we have the inequalities $$[M,M_{rs}]=\sum_{i=1}^{p}\sum_{j=1}^{n_{i}}[M_{ij},M_{rs}]\leq \sum _{i=1}^{p} l_{E_{i}}(M_{rs})n_{i}\leq \sum_{i=1}^{p}l_{E_{i}}(M_{rs})[U,E_{i}]=[U,M_{rs}].$$
Therefore the codimension is one iff we have 
$$[U,E_{i}]=n_{i} \mbox{  for all i with } l_{E_{i}}(M)\neq 0 \mbox{ and }$$$$[M_{ij},M_{rs}]= min \{l_{E_{i}}(M_{rs}),l_{Soc (M_{rs})}(M_{ij})\} = l_{E_{i}}(M_{rs})$$ for all $i,j,r,s$. This implies part a) immediately.

The assumption $\partial V \geq p$ implies $[U,E_{i}]>0$ and $[E_{i},V]>0$ for all $i$ by lemma \ref{gleichverteilung}. If all $E_{i}$ are composition factors and if the codimension is one, we have $[U,E_{i}]=n_{i}$ for all $i$. Thus all $E_{i}$ occur already in the regular top of $M$ and dually also in the regular socle. Furthermore, $[M_{ij},M_{rs}]$ is equal to $l_{E_{i}}(M_{rs})$ for all $i$, whence all $M_{ij}$'s with the same top $E_{i}$ have the same length. To see that the lengths also coincide for different regular tops, start with an index $i_{0}$ where this length $l$ is maximal and let $E_{j}$ be the regular socle of $M_{i_{0}k}$. Then $E_{j+1}$ occurs in the regular socle of some $M_{rs}$. It follows that $l(M_{rs}) \leq l$ by maximality. 
If $l(M_{rs}) < l-1$, we could proceed as in the proof of lemma \ref{gleichlang2} to derive a contradiction to the minimality of $M<U\oplus V$. The case $l(M_{rs}) =l-1$ is impossible since it would imply $Top(M_{rs})=E_{i_0}$.
Hence, we have $r=i_{0}+1$. So all $M_{i_{0}+1,k}$ also have maximal length. Our claim follows by induction. Here the indices of the simples $E_{i}$ have to be read modulo $p$. $\hfill\square$\\

\subsection{$M$ has less than $\partial V$ summands}

\begin{proposition}\label{sklein}
 Let $M<U\oplus V$ be a minimal degeneration such that $M$ belongs to a tube of period $p$ and has $s<\partial  V$ direct summands.  Then we have $Codim(U\oplus V,M)=1$.
\end{proposition}
\proof  First we look at the case $\partial V \leq p+1$. Suppose the codimension is not $1$. By lemma \ref{codimzugross} there is some  minimal degeneration $M'<U'\oplus V'$ and an indecomposable summand $X$ of $M$ occurring twice in $M'$ such that $[X,M]<[X,U\oplus V]$ and $$M\oplus X =M'\oplus Z < U'\oplus V' \oplus Z < U\oplus V \oplus X.$$  By lemma \ref{codimreg} $U'\oplus V'$ is  not regular. So we can assume that $U'$ is preprojective and that $V'$ is preinjective. For a regular simple $H$ with $\dimv(H)=\delta$ we have $\partial( V')=[H,V']\leq [H,V]=\partial(V)$ . Let $E$ be the regular socle of $X$. If $\partial(V')\leq p$, lemma \ref{gleichverteilung} delivers the contradiction $2\leq [E,X^{2}]\leq [E,M']\leq [E,V']\leq 1$.  So we can  assume that $\partial(V')= \partial(V)= p+1$. 

By assumption we have $s=\sum_{j=1}^{p}[M,E_{j}] < - \partial U= \sum_{j=1}^{p} [U,E_{j}]$. So we can fix an index $i$ with $[U,E_{i}]>[M,E_{i}]$. By lemma \ref{regkurze}, $E_{i}$ occurs only in the regular top of $M$. Thus $M$ is regular semi-simple for $p=1$ and the codimension is $1$ by proposition \ref{indanfang}. So we have $p\geq 2$ from now on.

Next, if $[M,E_{j}]=0$ or $[E_{j},M]=0$ for some $j$, then the codimension is one by lemma \ref{regkurze}. Hence because of $\partial V = p+1$  we can assume  $[M,E_{j}]=1=[E_{j},M]$ for all $j$.  Then each $E_{j}$ occurs exactly once as a regular top and once as a regular socle of a direct summand of $M$. Since $E_{i}$ occurs only in the regular top of $M$, the direct summand with regular top $E_{i+1}$ must be the simple $E_{i+1}$. Next, one looks at the indecomposable summand $W$ with regular top $E_{i+2}$. Since $E_{i}$ only occurs in the regular top of $M$ and $E_{i+1}$ occurs already with multiplicity $1$ in the regular socle, $W$ is $E_{i+2}$. Going on like that, one sees that $M$ is regular semi-simple and the codimension is $1$.

We are left with the case $\partial V\geq p+2$ and $[M,E_{j}]\neq 0 \neq [E_{j},M]$ for all $j$. First take $p=2$. Renumbering $i=1$ we have  $M=E_{1}^{a}\oplus E_{2}(2)^{b}\oplus E_{2}^{c}$. Assume $b>0$. Then $E_{1}$ occurs not only in the regular socle of $M$ and $E_{2}$ not only in the regular top. Consequently, lemma \ref{regkurze} implies $[E_{1},V]=[E_{1},M]=a$ and $c=[M,E_{2}]=[U,E_{2}]\geq 2$. It follows that $3\leq b+c=[E_{2},M]\leq [E_{2}V]\leq 3$. This gives the contradiction 
$$s=a+b+c=[E_{1},V]+[E_{2},V]=\partial V.$$ 
So $M$ is regular semi-simple and the codimension is $1$.

Finally, $p=3$ and $\partial V$ is $5$ or $6$. For $\partial V=6$ we have $[U,E_{j}]=2=[E_{j},U]$ for all $j$. Assume $i=1$ so that $E_{1}$ occurs only in the regular top and only once. Since $E_{1}$ occurs also in the regular socle, we have $M=E_{1}\oplus E_{2}^{a}\oplus E_{2}(2)^{b} \oplus E_{3}^{c}$. Here we have $a+b=[E_{2},M]\leq 2$. Therefore the degeneration $M'=E_{1}\oplus E_{2}^{a+b} \oplus E_{3}^{b+c}$ satisfies $M \leq M' \leq U \oplus V$. By minimality, $M=M'$ is regular semi-simple.

For $\partial V =5$ the situation is no longer symmetric and we use the notations from the table \ref{regsimplestable}. We have an $\tilde{E}_{8}$-quiver with $U=P(4)$ and $V=\tau^{k}I(4)$ for some $k$.
  From $s<\partial V$ we obtain $[U,S_{i}']=2 > [M,S_{i}']$ for $i=1$ or $3$. For $i=3$ one sees as above that $M$ is regular semi-simple. If $i=1$, $M$ must be regular semi-simple or isomorphic to $S_{1}'\oplus S_{2}'\oplus S_{2}'(2)\oplus S_{3}'$. But in the second case we have $q(\dimv(M)-\dimv(U))=2$. Hence, this case is excluded since there exists no preinjective indecomposable $\tau^{k}I(4)$ such that
$\dimv(U\oplus \tau^{k}V)= \dimv(M)$.
$\hfill\square$

\subsection{Some universal and generic extensions}

To finish the proof about the codimension we have to know some generic extensions that are given as universal extensions. To recall this concept, let $X,Y_{1},Y_{2}, \ldots Y_{r}$ be arbitrary modules with $n_{i}=[X,Y_{i}]^{1}$. Define $Y=\bigoplus_{i=1}^{r} Y_{i}^{n_{i}}=\bigoplus_{i=1}^{r} (\bigoplus_{j=1}^{n_{i}} Y_{ij})$ with $Y_{ij}=Y_{i}$ for all $i$ and $j$. Now $\epsilon$ in $Ext(X,Y)$ is called universal if the components $\epsilon_{i1},\epsilon_{i2},\ldots ,\epsilon_{in_{i}}$ are a  basis of $Ext(X,Y_{i})$ for all $i$.
\begin{lemma}\label{uniext}
 Under the above assumptions we have:
\begin{enumerate}
 \item All universal extensions are conjugate under the natural action of $Aut(Y)$ on $Ext(X,Y)$. In particular, all middle terms of universal extensions are isomorphic to a fixed module $Z$. This is the generic extension.
\item If $X=\bigoplus_{j=1}^{m}X_{j}$, then the generic extension $Z_{j}$ of $X_{j}$ by $\bigoplus_{i=1}^{r}Y_{i}^{[X_{j},Y_{i}]^{1}}$ exists for all $j$ and we have $Z=\bigoplus_{j=1}^{m} Z_{j}$.
\end{enumerate}
\end{lemma}
\proof a) The group $G=\prod GL_{n_{i}}$ can be diagonally embedded into $Aut(Y)$. It acts transitively on the set of universal extensions. Any extension is a degeneration of $Z$ because the set of ordered bases $(b_{1},\ldots ,b_{t})$ in $k^{t}$ is dense in the set of all $t$-tuples of vectors.

b) The direct sum of universal extensions $\epsilon_{j}$ of $X_{j}$ by $\bigoplus_{i=1}^{r}Y_{i}^{[X_{j},Y_{i}]^{1}}$ is a universal extension of $X$ by $Y$, whence their middle terms are isomorphic.$\hfill \square$\\

We are interested in the special case where $\mathcal{T}$ is a tube with regular simples $E_{1},E_{2},\ldots E_{p}$. For any module $X$ we denote by $e(X)$ the middle term of the universal extension of $X$ by $\bigoplus _{i=1}^{p}E_{i}^{[X,E_{i}]^{1}}$. By part b) of the last lemma we only have to know $e(X)$ for indecomposable $X$.

\begin{proposition}\label{genext} Under the above assumptions we have:
\begin{enumerate}
 \item For $X$ indecomposable in the given tube $e(X)$ is the indecomposable with the same regular top as $X$ and with regular length $l(X)+1$. One has $[X,E_{k}]^{1}=[E_{k},e(X)]=[e(X),E_{k-1}]^{1}$ for all $k$. Thus $e^{p}(X)$ is the indecomposable with the same top as $X$ and dimension vector $\dimv(X)+\delta$.
\item If $X$ indecomposable preinjective, $e(X)$ is the preinjective indecomposable with dimension vector $\dimv(X) + \sum_{i}[X,E_{i}]^{1}\dimv(E_{i})$. Again, we have $[X,E_{k}]^{1}=[E_{k},e(X)]^1=[e(X),E_{k-1}]$ for all $k$. Thus $e^{\epsilon(Q)p}(X)$ is $\tau^{p(Q)}(X)$.
\item For all other indecomposables $e(X)=X$ is true.
\item If $X$ is preinjective indecomposable with $dim \,X > dim \bigoplus_{i=1}^{p}E_{i}^{[E_{i},X]}$, there is an exact sequence $0 \longrightarrow \bigoplus_{i=1}^{p}E_{i}^{[E_{i},X]} \longrightarrow X \longrightarrow X' \longrightarrow 0$ with an indecomposable preinjective $X'$. Here we have $e(X')=X$.
\end{enumerate}
\end{proposition}
\proof Only b) and d) are  not obvious.

b) We set $E=\bigoplus _{i=1}^{p}E_{i}^{[X,E_{i}]^{1}}$ and we calculate $[X,X]-[X,X]^{1}+[E,X]-[X,E]^{1}+[E,E]-[E,E]^{1}=1$ because $[E,E]=[X,E]^{1}=\sum_{i}([X,E_{i}]^{1})^{2}$ and $[E,X]=[E,E]^{1}=[\tau^{-}E,E]=\sum_{i}[E_{i},X][E_{i+1},X].$ Thus $\dimv (X) + \dimv(E)$ is the dimension vector of a preinjective indecomposable $Z'$.
 Let $0 \longrightarrow \bigoplus_{i=1}^{p}E_{i}^{[X,E_{i}]^{1}} \longrightarrow Z \longrightarrow X \longrightarrow 0$ be the universal extension.  From $[Z,F]\leq [E,F]$ for all regular simples we get that $Z_{R}$ belongs to the tube $\mathcal{T}$. In the induced exact sequence 
$$0\longrightarrow Hom(Z,E_{k}) \rightarrow Hom (E,E_{k}) \rightarrow Ext(X,E_{k})$$ the last map is  an isomorphism for all $E_{k}$ by the definition of the universal extension. Thus also $[Z,E_{k}]=0$ for all $E_{k}$ whence $Z$ is preinjective. 

Clearly we have $[E,Z']=[E,Z]$ and  $Z' \leq Z$. Since the orbit of 
$X$ is open, $Z'$ is an extension  by \cite[lemma 4.4]{Degenerations}. Since $Z$ is the generic extension, we have $Z=Z'$.

We calculate $[e(X),E_{k-1}]^{1}=[e(X),E_{k-1}]^{1}-[e(X),E_{k-1}]=[X,E_{k-1}]^{1}-[X,E_{k-1}]+[E,E_{k-1}]^{1}-[E,E_{k-1}]=[X,E_{k}]^{1}=:n_{k}$. This implies that $e^{\epsilon(Q)p}(X)$ has dimension vector $\dimv(X) + (\epsilon(Q)\sum_{i=1}^{p}n_{i})\delta$. Since $\sum_{i=1}^{p}n_{i}=\sum_{i=1}^{p}[E_{i+1},X])=\partial(X)$, this is also the dimension vector of $\tau^{p(Q)}(X)$. Thus $\tau^{p(Q)}(X)=e^{\epsilon(Q)p}(X)$ since both are indecomposable preinjective.

d) As in b), one checks that $\dimv X -\dimv \bigoplus_{i=1}^{p}E_{i}^{[E_{i},X]}$ is the dimension vector of a preinjective indecomposable $X'$. We have $[X',E_{k}]^{1}=[X',E_{k}]^{1}-[X',E_{k}]=[X,E_{k}]^{1}-[X,E_{k}]-[\bigoplus_{i=1}^{p}E_{i}^{[E_{i},X]},E_{k}]^{1}+[\bigoplus_{i=1}^{p}E_{i}^{[E_{i},X]},E_{k}]=[E_{k},X]$. The claim follows from part b).$\hfill \square$\\

\subsection{$M$ has $\partial V$ direct summands}

\begin{proposition}\label{indsub}
 Let $M<U\oplus V$ be a minimal degeneration such that $M$ is the direct sum of $\partial V$ indecomposable summands. Suppose $M$ is not regular semi-simple. Let $E$ be the regular socle of $M$. Set $Q=M/E$ and $M'=U\oplus V$. Then the generic quotient $V'$ of $V$ by $E$ is indecomposable preinjective and $Q<U\oplus V'$ is a minimal degeneration of the same codimension as $M<U\oplus V$.
\end{proposition}
\proof We want to apply theorem \ref{divsub} and check that the assumptions are all satisfied. Since $M$ has $\partial V$ indecomposable summands, we have $[E_{i},M]=[E_{i},V]$ for all $i$ and therefore $E=\bigoplus_{i=1}^{p}E_{i}^{[E_{i},V]}$. This implies condition i). Clearly, $Q$ is the only quotient of $M$ by $E$ and $M$ is an extension of $Q$ by $E$ with $M\leq M'$. If $M''$ is another extension with $M'' \leq M'$ one has $[E_{j},V] \leq [E_{j},M'']\leq [E_{j},M']=[E_{j},V]$ for all regular simples $E_{j}$. Thus $M''$ has the same regular socle $E$ as $M$ and the same quotient. Therefore $M$ and $M''$ are isomorphic. So $M$ is the generic extension and the only one degenerating to $M'$. 

Since $M$ is not regular semi-simple we have $dim V > dim E$. By lemma \ref{genext} there is an exact sequence $0 \rightarrow E \rightarrow V \rightarrow V' \rightarrow 0$ with indecomposable preinjective $V'$. This implies condition iii). 

 Finally the assumptions about generic extensions in iv) and generic quotients in v) hold because there are always only finitely many isomorphism classes as candidates.$\hfill \square$

\begin{proposition}
 If $M<U\oplus V$ is a minimal degeneration with $M$ in one tube, one has $Codim(U\oplus V,M)=1$.
\end{proposition}
\proof We proceed  by induction on $dim \,M$. If $M$ is regular semi-simple or if $s<\partial V$, the assertion is true by proposition \ref{indanfang} or by proposition \ref{sklein}. If $s= \partial V$, we simply apply proposition \ref{indsub} and the induction hypothesis.$\hfill \square$\\

This completes the proof of theorem \ref{Theo1}.

\section{On the classification} \label{classsec}

\subsection{Degenerations between preprojective modules}

In \cite{Fritzsche} there is shown a periodic behaviour of the minimal degenerations that reduces their classification to a finite problem. This is then solved by using a computer. In particular, the indigestible lists show that the codimensions of the building blocs are bounded by $2$. This is proved  now in proposition \ref{codimpre} by theoretical means.

\subsection{$U$ is preprojective and $V$ regular}

By part a) of proposition \ref{codimVreg}, the classification of the minimal disjoint degenerations is reduced to the case where $\dimv(V)\leq \delta$. This case is treated  using a computer in \cite{Wolters}. To get an impression of the arising complexity  one can look at section \ref{defext}.

\subsection{ $U$ is preprojective and $V$ preinjective}
In the special case where $-\partial(U)=\partial(V)=1$ the degeneration test \ref{enttest} implies the classification of the building blocs.  Note that this settles the case $Q=\tilde{A}_m$.

\begin{proposition}\label{defectone} Supposed that $-\partial(U)=\partial(V)=1$,
a module $M$  degenerates to $U\oplus V$ if and only if $\dimv(M)=\dimv(U\oplus V)$ and $M$ is the direct sum of indecomposables $M_i$ from pairwise different tubes such that $[U,Top(M_i)]=1$.
In particular, any degeneration is minimal.
\end{proposition}

If $-\partial(U)\geq 2$ or $\partial(V)\geq 2$, the situation is much more complicated. Anyhow,
by exploiting the technique of shrinking or extending suitable arrows, we may assume that $Q$ is of type $\tilde{D}_8$, $\tilde{E}_6$, $\tilde{E}_7$ or $\tilde{E}_8$, and thus, that we have $\epsilon(Q)=1$.

Provided that $M$ has enough direct summands in one tube,  we can point out a periodic behaviour for the building blocs.
\begin{proposition}\label{periprop}
Let  $Q$ and $Q'=U\oplus V$ be modules with the same dimension vector such that $[E_{i},Q]=[E_{i},V]$ holds for all regular simples in the tube $\T$. Then we have:
\begin{enumerate}
\item $Q\leq Q'$ iff $e(Q)\leq e(Q')$.
\item $Codim(e(Q'),e(Q))=Codim(Q',Q)$.
 \item $Q<Q'$ is minimal iff $e(Q)<e(Q')$ is minimal.
\item $[E_{i},e(Q)]=[E_{i},e(Q')]$ for all $i$.
\end{enumerate}
\end{proposition}
\proof  We want to apply theorem \ref{divsub} to $E=\bigoplus_{i=1}^{p}E_{i}^{[E_{i},e(Q)]}$, $M=e(Q)$, $M'=e(Q')$, $Q$ and $Q'$. Recall that $n_i:=[E_{i},e(Q)]=[Q,E_{i}]^{1}$ always holds. Proposition \ref{genext} implies $[E_{i},e(Q)]=[E_{i},e(Q')]$ for all $i$, whence part d) above and  $[E,M]=[E,M']$ are true. By construction, $E$ embeds into $M$ resp. $M'$. Let $\epsilon$  be such an embedding of $E$ into $M$ resp. $M'$. Then the components $\epsilon_{1,1}, \ldots , \epsilon_{1,n_1}, \epsilon_{2,1}, \ldots , \epsilon_{p,n_p}$ of $\epsilon$ satisfy 
\begin{eqnarray*}
 dim(E)&=&dim(Im(\epsilon))=dim(\sum_{i,j}Im(\epsilon_{i,j}))\\&\leq &\sum_{i,j}dim(Im(\epsilon_{i,j}))
\leq \sum_{i}n_idim(E_i)=dim(E),
\end{eqnarray*}
whence the $\epsilon_{i,1}, \ldots \epsilon_{i,n_i}$ are linearly independent for all $i$. Since $n_i=[E_i,M]=[E_i,M']$, these maps actually form a basis of $Hom(E_i,M)$ resp. $Hom(E_i,M')$. Hence, the image of an embedding $\epsilon$ of $E$ into $M$ resp. $M'$ is in both cases equal to the sum of the images of all homomorphisms from  $E$ to $M$ resp. $M'$. Consequently, there is in both cases only one quotient. Due to proposition \ref{genext}, $M$ and $M'$ are the generic extensions.
So parts a) and b) follow from theorem \ref{divsub}. Since we know already that the codimensions of all building blocs are $1$, part c) is also true.
But the remainig conditions of theorem \ref{divsub} are also satisfied as one can show.\nolinebreak $\hfill \square$

\begin{corollary}{\bf (Periodicity Theorem)}\label{peritheo}
Let  $M_P$ be preprojective, $M_R$ be regular with no summand belonging to the tube $\T$ and $M_I$ be preinjective such that 
\[ \partial(M_P)-\partial(U)=\partial(V)-\partial(M_I)>0.\]
If $R \in\T$ satisfies $[U,E_i]=[M_P\oplus R,E_i]$ for all simples $E_i\in \T$, 
 then the following statements are equivalent.
\begin{enumerate}[(i)]
\item $M_P\oplus M_R\oplus R\oplus M_I<U\oplus V$ is minimal. 
\item $M_P\oplus M_R\oplus e^{p}(R)\oplus \tau^{p(Q)}M_I<U\oplus \tau^{p(Q)}V$ is minimal.
 
\end{enumerate}
\end{corollary}
\proof This follows  from proposition \ref{periprop} and proposition \ref{genext}. $\hfill\square$

Recall, the easy description of $e^p(R)$ is given in part a) of proposition \ref{genext}.
Of course, supposed that $\T$ is homogeneous with regular simple $E$ and that one of the modules $M_P$, $M_I$ or $M_R$ is non-zero, the assertion of the periodicity theorem also holds if one replaces  $R$ by $0$ and $e^{p}(R)$ by $E^{\partial(V)-\partial(M_I)}$.

For the preprojective resp. the preinjective part of 
$M$, one obtains the following lemma by using essentially the same technical arguments as in  \cite[lemma 4]{Fritzsche}. The full proof is contained in \cite{Frank}.

\begin{lemma} \label{naheprepsum} If $M<U\oplus V$ is minimal, then it holds:
\begin{enumerate}
 \item For every direct summand $X$ of $M_P$ it holds
$
 d(U,X)< 4p(Q)+d(Q).
$
\item For every direct summand $X$ of $M_I$ we have
$
 d(X,V)<
4p(Q)+d(Q).
$
\end{enumerate}
\end{lemma}

\begin{corollary} 
 Up to application of the periodicity theorem \ref{peritheo} the number of minimal degenerations to the direct sum of a simple projective and  preinjective indecomposable  is finite.
\end{corollary}
\proof 
Let $M'<U\oplus V'$ be minimal such that $ t\partial(V')\delta \leq \dimv(V')< (t+1)\partial(V')\delta$ for some $t\geq 78$. To prove the assertion, it suffices to show the existence of some tube $\T$ such that $R':=M_\mu'$ has $\partial(V')-\partial(M_I')$ indecomposable summands $R_i'$ of length $>p_\mu$ (resp. $\geq p_\mu$ for homogeneous $\T$). 

Assuming the contrary, part c) of lemma \ref{gleichlang} insures that there is no summand coming from a homogeneous tube. Therefore, we can write
\[ M'=M_P'\oplus M_1'\oplus M_2' \oplus M_3'\oplus M_I'\] where  the $M_i'$ are the summands from the three non-homogeneous tubes $\mathcal{T}_1$, $\mathcal{T}_2$, and $\mathcal{T}_3$. For any indecomposable direct summand $X$ of $M_i'$  lemma \ref{gleichlang2} resp. lemma \ref{regkurze} yield $l(X)\leq 2p_i$.
Due to lemma \ref{naheprepsum} any summand $X$ of $M_P'$ satisfies $d(U,X)< 4p(Q)+d(Q)$, whence we can generously estimate  $\dimv(M_P')\leq -4\partial(M_P')\delta\leq 20\delta$. The same lemma implies $\dimv(M_I')\leq (t+4)\partial(M_I')\delta\leq (t+4)(\partial(V')-1)\delta$.
Therefore, we obtain
\begin{eqnarray*}
 0&=& \dimv(U\oplus V')-\dimv(M')\\
&\geq & t\partial(V')\delta - 56\delta - (t+4)(\partial(V')-1)\delta\\
&= &t\delta - 76\delta  \geq 2 \delta,
\end{eqnarray*}
which is a contradiction. Hence there exists some minimal degeneration $M:=M_P\oplus M_R\oplus R\oplus M_I<U\oplus V$ with $M_P'=M_P$, $M_R'=M_R$, $R'=e^p(R)$, $M_I'=\tau^{p(Q)}M_I$ resp. $V'=\tau^{p(Q)}V$ and $M'<U\oplus V'$ is obtained from $M<U\oplus V$ by exploiting the periodicity. $\hfill\square$ \\

This finishes the proof of theorem \ref{Theo2}.



\subsubsection{M in one tube}\label{classtubesec}

If $M$ comes from a single non-homogeneous tube $\mathcal{T}$, the building blocs are classified by means of a computer. Up to application of the propositions \ref{indsub} and \ref{periprop} these degenerations are listed in section \ref{mindeglist}. For a more detailed list, see \cite{Frank}. 
Besides, we make several interesting observations.
\begin{enumerate}[i)]
 \item In fact, there occur no building blocs $Q<U\oplus V'$ such that $Q\in \mathcal{T}$ has less than $\partial(V')$ direct summands, except those that arise from application of proposition \ref{indsub} to some $M<U\oplus V$ where $M$ has the required $\partial(V)$ indecomposable summands. For example, if $Q$ is of type $\tilde{E}_8$, the building bloc $S_4<P(9)\oplus \tau I(9)$ is produced in this way from $S_2\oplus S_3(2)\oplus S_5< P(9)\oplus \tau^{7}I(9)$.
\item The set $\mathcal{D}:=\{M\in \mathcal{T} \vert M<U\oplus V\}$ is either empty or it contains exactly one maximal element $M_0$. 
If $M_0<U\oplus V$ is not a minimal degeneration, then we have $\partial(V)\geq p+1$.
\item The propositions \ref{indsub} and \ref{periprop} have further refinements. For instance, proposition \ref{indsub} still holds if one replaces $E$ by a submodule $E'=\bigoplus_{j\in J}E_j^{n_j}$ of $E$ such that $q(\dimv(V)-\dimv(E'))=1$ and $n_j=[E_j,M]=[E_j,V]$ for all $j\in J\subseteq \{1,\ldots , p\}$. 
With regard to that,  for a fixed simple projective $U$ there exists essentially one building bloc in each tube. As an example, we choose $Q$ of type $\tilde{E}_6$,  $U=P(2)$, $V=\tau^2I(2)$, $M=S_1(1)\oplus S_2(2)$ and $E'=S_2$. Then $V'=V/E'=\tau^2 I(6)$ and $ S_1(1)\oplus S_3(1)  <U\oplus V'$ is again minimal.
\end{enumerate}

\section{Concluding remarks} \label{concludingsec}

\subsection{Tame concealed algebras}
Our main result about the codimension holds also for tame concealed algebras as we will explain now. The definition and the basic properties of tame concealed algebras  can be found in Ringels nice book \cite{Ringel}.
Here it is  important that by \cite{Degenerations,Extended,Extensions} for modules over a tame concealed algebra $B$ any minimal disjoint degeneration stems from a building bloc $M<U\oplus V$ as before. We claim that the codimension is bounded by two. This is clear by proposition  \ref{codimreg} if $U$ and $V$ are regular.

  We look at the slice in the preprojective component with $U$ as its only source. If this slice hits all $\tau$- orbits, we take the corresponding tilting module $T$ and its endomorphism algebra  which is a tame  path algebra. Since all indecomposables occurring in $M$ and $U\oplus V$ are generated by $T$, our problem is transferred to the known situation of extended Dynkin-quivers. If the slice does not hit all $\tau$-orbits, we can complete the partial tilting module given by the slice with some projectives to a tilting module $T$. For the new tame concealed algebra $B$ we have reduced to the case where $U$ is simple projective. 

If $V$ is also preprojective or regular, both modules are images under the functor $F=Hom (T,-)$ of indecomposable torsion modules $U',V'$ over a tame path algebra $A$. Then any deformation $M$ is  of the form $M=FM'$ for some deformation $M'$ of $U'\oplus V'$ which is also a torsion module. The codimensions coincide by the tilting theorem.

Next suppose that $V$ is preinjective. If the support of $V$ is $B$, we can dualize and look at the slice in the preinjective component with $V$ as the only sink. Taking the corresponding tilting module and its endomorphism algebra we are back in a quiver situation. So let the support of $V$ be a proper subalgebra $C$. If the support of $U$ is not contained in $C$, the codimension is  $1$ because $[U,V]=0$. If the support of $U$ lies in $C$, we are dealing with a representation directed algebra having only Dynkin-quivers as sections in the Auslander-Reiten quiver. The codimension is again $1$ as one sees by generalizing \ref{redtheo} and \ref{decofcodim} to the present situation.

\subsection{Minimal singularities}
In \cite{codim2} Zwara has classified the codimension two singularities occurring for representations of tame quivers. It follows from this article that his result covers all minimal disjoint singularities. So it is a natural question to ask which general minimal singularities occur.

\subsection{Which deformations are extensions}\label{defext}

Given two modules $U$ and $V$, it is an interesting open problem to determine which modules $M$ appear in the middle of a short exact sequence $0\ \rightarrow U \rightarrow M \rightarrow V \rightarrow 0$. The obvious necessary condition that $M$ is a deformation of $U\oplus V$ is in general not sufficient. We will discuss here the case where $U$ and $V$ are indecomposable representations of a tame quiver. 

If none of the modules is regular, any deformation is an extension by \cite[theorem 4.5]{Degenerations}. For regular $U$ and $V$ the situation is completely analyzed in \cite{Wolters}. Roughly speaking, only half of the deformations are extensions. 

\begin{figure}[hb!]
\centering
\includegraphics[width=0.6\textwidth]
{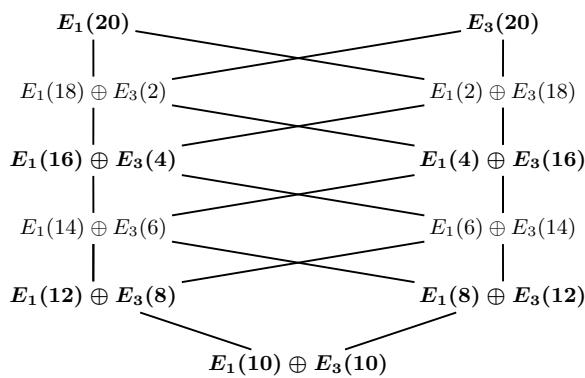}
\caption{Hasse diagram of deformations of $E_{1}(10)\oplus E_{3}(10)$ in a tube of period $4$ }\label{reghasse}
\end{figure}

As an example we take a tube with $4$ simples and $U=E_{1}(10)$, $V=E_{3}(10)$. Figure \ref{reghasse} illustrates the Hasse diagram of deformations of $E_{1}(10)\oplus E_{3}(10)$. In the diagram
only the bold deformations are extensions of $V$ by $U$ or $U$ by $V$ as follows from lemma \ref{codimreg}. In particular, the set of all extensions is far from being locally closed.

Next we take $U$ simple projective and $V$ regular with $\dimv(V)\leq \delta$. If $V$ is homogeneous or  $\dimv(V)<  \delta$, any deformation is again an extension as follows from the dual of \cite[theorem 2.4]{Degenerations}.

\begin{figure}[ht!]
\includegraphics[width=1.0\textwidth]{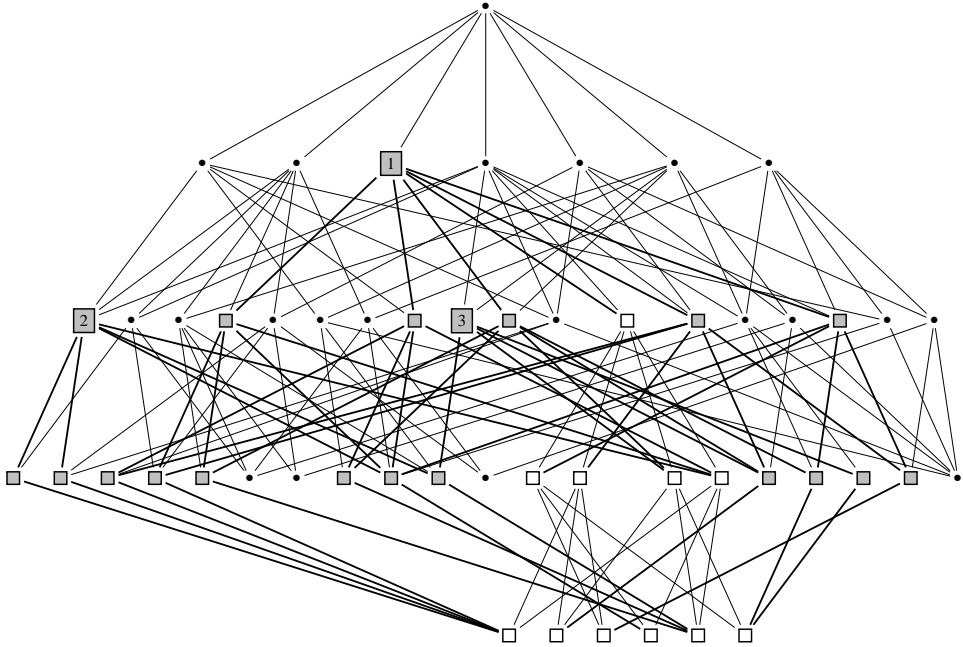}
\caption{Hasse diagram of proper deformations of $P(4) \oplus S_2(4)$}
\label{diaIsabel}
\end{figure}

\begin{proposition}
 Let $V$ be a regular indecomposable with $\dimv(V)= \delta$ in a tube with $p\geq 2$. Define $K$ to be the set of all deformations $M$ of $U \oplus V$ such that $M_{R}\neq 0$ or $M\leq M'$ where $M'$ is a minimal deformation of $U\oplus V$ which is preprojective. Then $K$ is the  set of all extensions.
\end{proposition}
\proof Using again results from \cite{Degenerations} one gets quite easily that $K$ consists only of extensions. The other inclusion is shown in \cite{Wolters} case by case using a computer as well as some handwork. The result is illustrated  in figure \ref{diaIsabel}.$\hfill \square$\\

Figure \ref{diaIsabel} shows the Hasse-diagram of the proper deformations of $P(4)\oplus S_{3}(4)$ for the quiver $\tilde{E}_{7}$ ( see table \ref{regsimplestable} ). The codimensions are 1,2,3,4,5 from the bottom to the top where the indecomposable $\tau^{-3}P(4)$ is. The white empty boxes represent deformations $M$ with $M_{R}\neq 0$. The black points are preprojective deformations $M$ such that there is a minimal deformation $M'$ with $M\leq M'$ and $M'=M'_{P}$. These two types of deformations are the only extensions. The large grey boxes represent the modules $M_{1}=\tau^{-3}P(2)\oplus \tau^{-1}P(2), M_{2}=\tau^{-4}P(2)\oplus P(6), M_{3}=\tau^{-4}P(7)\oplus \tau^{-2}P(8)\oplus P(1)$. These are the maximal deformations that are not extensions as one checks by hand.

\subsection{Codimensions and representation type}

Given any algebra $A$, let us denote by $Cod(A)$ the supremum of the codimensions of all minimal disjoint degenerations $M\leq N$ of $A$-modules. Then we have:
\begin{proposition} For a path algebra $A$ of a connected quiver $Q$ it follows:
\begin{enumerate}
 \item $Cod(A)=1$ iff $Q$ is a Dynkin-quiver.
\item $Cod(A)=2$ iff $Q$ is Euclidean.
\item $Cod(A)=\infty$ iff $Q$ is wild.
\end{enumerate}
\end{proposition}
\proof The first two parts follow from our results. For the generalized Kronecker quiver with three arrows one gets easily minimal disjoint degenerations between preprojectives of arbitrarily high codimension. The last part is then clear, since the embedding of any module category into the module category of a wild quiver can be done in a way compatible with degenerations by \cite[example 5.19]{trond}.$\hfill \square$\\

In \cite{Olbricht} there are examples of wild quivers such that for any natural number $n\geq 2$ there are preprojective indecomposables $U$ and $V$ and minimal deformations $M_{i}$ satisfying $Codim(U\oplus V,M_{i})=i$ for all $i$ between $2$ and $n$. Thus the degenerations of modules over wild quivers behave as they should: wild!

 \section{Tables}\label{tablessec}

In this last section, we append the tables that we referred to during this article. 
Let $Q$ be a quiver of type $\tilde{D}_8$, $\tilde{E}_6$, $\tilde{E}_7$ or $\tilde{E}_8$. Recall, by tilting we reduced to the case where $Q$ has only one sink $i$, i.e. there is exactly one simple projective $U:=P(i)$. Thus $Q$ is uniquely determined by the type of $Q$ together with the specification of $U$. The following lists make use of this fact.

\subsection{The quivers}

The points of $Q$ are numbered in the following way.
\vspace{0.4cm}
\begin{longtable}{rl}
 \endhead
\endfoot
\mbox{$ \tilde{D}_8:$}&\parbox{9.2cm}{$\xymatrix{ 1&&&&&&8\\
&3\ar@{-}[lu]\ar@{-}[r]\ar@{-}[dl]& 4 \ar@{-}[r]   & 5\ar@{-}[r] &  6 \ar@{-}[r]&7\;\;\ar@{-}[l]\ar@{-}[ur]\ar@{-}[dr]&\\
2&&&&&&9\\
}$}\\\vspace{0.2cm}
$\tilde{E}_6$:&\parbox{9.2cm}{$\xymatrix{&&7\ar@{-}[d]&&\\
 &&6\ar@{-}[d]&&\\
1&2\ar@{-}[l]\ar@{-}[r]&3\ar@{-}[r]&4\ar@{-}[r]&5
}$}\\\vspace{0.2cm}
$\tilde{E}_7$:&\parbox{9.2cm}{$\xymatrix{
 &&&8\ar@{-}[d]&&\\
1&2\ar@{-}[l]\ar@{-}[r]&3\ar@{-}[r]&4\ar@{-}[r]&5\ar@{-}[r]&6\ar@{-}[r]&7}
$}\\\vspace{0.2cm}
$\tilde{E}_8$:&\parbox{9.2cm}{$\xymatrix{
 &&9\ar@{-}[d]&&\\
1&2\ar@{-}[l]\ar@{-}[r]&3\ar@{-}[r]&4\ar@{-}[r]&5\ar@{-}[r]&6\ar@{-}[r]&7\ar@{-}[r]&8}
$}\\
\end{longtable}

\input{regsimples}

\input{tables}

\end{document}

%% file: regsimples.tex
\subsection{The regular simples}\label{regsimplestable}

The next table gives dimension vectors of the regular simples of the three non-homogeneous tubes $\mathcal{T}_k$, $1\leq k\leq 3$. 
Here, the simples of the tube $\mathcal{T}_1$, $\mathcal{T}_2$ resp. $\mathcal{T}_3$ are denoted by $S_i$, $S_i'$ resp. $S_i''$.
As before, the convention is $\tau S_i=S_{i-1}$. 
\renewcommand{\arraystretch}{1.6}
\begin{longtable}{||c|c|c|c|c|c|c||}\hline\hline
$\!\!|Q|\!$&$U$&$k$&$\!p_{\,\!k}\!\!\!  $&\multicolumn{3}{l||}{\mbox{\small Regular simples}}\\\hline
\endhead
\hline\hline
\endfoot
%
%
$\!\!\!\tilde{D}_8\!\!\!$&$\!\!P(3)\!\!$&$1$&$6$&$\!{    S_1\!=\!\lle \arr{0}{0}\!\!  0 \!\; 0 \!\; 0 \!\; 1 \!\; 0 \!\!\arr{0}{0}\rri   }$&$ \!{    S_2\!=\!\lle \arr{0}{0}\!\! 0 \!\; 0 \!\; 0 \!\; 0 \!\; 1  \!\!\arr{0}{0}\rri   }$&$ {     {  \!S_3\!=\!\lle \arr{0}{0}\!\! 1 \!\; 1 \!\; 1 \!\; 1 \!\; 1  \!\!\arr{1}{1}\rri   }}$ \\
&&&&${     {  \!S_4\!=\!\lle \arr{1}{1}\!\! 1 \!\; 0 \!\; 0 \!\; 0 \!\; 0  \!\!\arr{0}{0}\rri   }}$&$ \!{    S_5\!=\!\lle \arr{0}{0}\!\! 0 \!\; 1 \!\; 0 \!\; 0 \!\; 0  \!\!\arr{0}{0}\rri   }$&$ \!{    S_6\!=\!\lle \arr{0}{0}\!\! 0 \!\; 0 \!\; 1 \!\; 0 \!\; 0  \!\!\arr{0}{0}\rri   }$  \\
&&$2$&$2$&${     {  \!S_1'\!=\!\lle \arr{1}{0}\!\! 1 \!\; 1 \!\; 1 \!\; 1 \!\; 1  \!\!\arr{1}{0}\rri   }}$&$ {     {  \!S_2'\!=\!\lle \arr{0}{1}\!\! 1 \!\; 1 \!\; 1 \!\; 1 \!\; 1  \!\!\arr{0}{1}\rri   }}$&  \\
 &&$3$&$2$&${     {  \!S_1''\!=\!\lle \arr{1}{0}\!\! 1 \!\; 1 \!\; 1 \!\; 1 \!\; 1  \!\!\arr{0}{1}\rri   }}$&$ {     {  \!S_2''\!=\!\lle \arr{0}{1}\!\! 1 \!\; 1 \!\; 1 \!\; 1 \!\; 1  \!\!\arr{1}{0}\rri   }}$&  \\\hline
&$\!\!P(4)\!\!$&$1$&$6$&$\!{    S_1\!=\!\lle \arr{0}{0}\!\! 0 \!\; 0 \!\; 0 \!\; 1 \!\; 0  \!\!\arr{0}{0}\rri   }$&$ \!{    S_2\!=\!\lle \arr{0}{0}\!\! 0 \!\; 0 \!\; 0 \!\; 0 \!\; 1  \!\!\arr{0}{0}\rri   }$&$ {     {  \! S_3\!=\!\lle \arr{0}{0}\!\! 0 \!\; 1 \!\; 1 \!\; 1 \!\; 1  \!\!\arr{1}{1}\rri   }}$\\
&&&&$ \!{    S_4\!=\!\lle \arr{0}{0}\!\! 1 \!\; 0 \!\; 0 \!\; 0 \!\; 0  \!\!\arr{0}{0}\rri   }$&$ {     {  \! S_5\!=\!\lle \arr{1}{1}\!\! 1 \!\; 1 \!\; 0 \!\; 0 \!\; 0  \!\!\arr{0}{0}\rri   }}$&$ \!{    S_6\!=\!\lle \arr{0}{0}\!\! 0 \!\; 0 \!\; 1 \!\; 0 \!\; 0  \!\!\arr{0}{0}\rri   }$ \\
&&$2$&$2$&${     {  \! S_1'\!=\!\lle \arr{1}{0}\!\! 1 \!\; 1 \!\; 1 \!\; 1 \!\; 1  \!\!\arr{1}{0}\rri   }}$&$ {     {  \!S_2'\!=\!\lle \arr{0}{1}\!\! 1 \!\; 1 \!\; 1 \!\; 1 \!\; 1  \!\!\arr{0}{1}\rri   }}$&  \\
&&$3$&$2$&${     {  \!S_1''\!=\!\lle \arr{1}{0}\!\! 1 \!\; 1 \!\; 1 \!\; 1 \!\; 1  \!\!\arr{1}{0}\rri   }}$&$ {     {  \! S_2''\!=\!\lle \arr{0}{1}\!\! 1 \!\; 1 \!\; 1 \!\; 1 \!\; 1  \!\!\arr{0}{1}\rri   }}$&  \\\hline
&$\!\!P(5)\!\!$&$1$&$6$&$\!{    S_1\!=\!\lle \arr{0}{0}\!\! 0 \!\; 0 \!\; 0 \!\; 1 \!\; 0  \!\!\arr{0}{0}\rri   }$&$ \!{    S_2\!=\!\lle \arr{0}{0}\!\! 0 \!\; 0 \!\; 0 \!\; 0 \!\; 1  \!\!\arr{0}{0}\rri   }$&$ {     {  \! S_3\!=\!\lle \arr{0}{0}\!\! 0 \!\; 0 \!\; 1 \!\; 1 \!\; 1  \!\!\arr{1}{1}\rri   }}$\\
&&&&$ \!{    S_4\!=\!\lle \arr{0}{0}\!\! 0 \!\; 1 \!\; 0 \!\; 0 \!\; 0  \!\!\arr{0}{0}\rri   }$&$ \!{    S_5\!=\!\lle \arr{0}{0}\!\! 1 \!\; 0 \!\; 0 \!\; 0 \!\; 0  \!\!\arr{0}{0}\rri   }$&$ {     {  \! S_6\!=\!\lle \arr{1}{1}\!\! 1 \!\; 1 \!\; 1 \!\; 0 \!\; 0  \!\!\arr{0}{0}\rri   }}$  \\
&&$2$&$2$&${     {  \!S_1'\!=\!\lle \arr{1}{0}\!\! 1 \!\; 1 \!\; 1 \!\; 1 \!\; 1  \!\!\arr{1}{0}\rri   }}$&$ {     {  \! S_2'\!=\!\lle \arr{0}{1}\!\! 1 \!\; 1 \!\; 1 \!\; 1 \!\; 1  \!\!\arr{0}{1}\rri   }}$&  \\
&&$3$&$2$&${     {  \! S_1''\!=\!\lle \arr{1}{0}\!\! 1 \!\; 1 \!\; 1 \!\; 1 \!\; 1  \!\!\arr{0}{1}\rri   }}$&$ {     {  \!S_2''\!=\!\lle \arr{0}{1}\!\! 1 \!\; 1 \!\; 1 \!\; 1 \!\; 1  \!\!\arr{1}{0}\rri   }}$&  \\\hline\hline
%
%
&&&&&&\vspace*{-0.615cm}\\
$\!\!\!\tilde{E}_6 \!\!\!$&$\!\!P(2)\!\!$&$1$&$ 3$&${     {  \!S_1\!=\!\lle {\arw\begin{array}{ccccc} &\moem&\moem0&\moem&\moem\\&\moem&\moem0&\moem&\moem\\ 1 &\moem 1 &\moem 1 &\moem 1 &\moem 1 \end{array} } \rri   }}$&$ \!{    S_2\!=\!\lle {\arw \begin{array}{ccccc} &\moem&\moem 0 &\moem&\moem \\ &\moem&\moem 1 &\moem&\moem\\  0 &\moem 0 &\moem 1 &\moem 0 &\moem 0     \end{array} } \rri   }$&$ {     {  \! S_3\!=\!\lle {\arw \begin{array}{ccccc}&\moem&\moem 1 &\moem&\moem \\ &\moem&\moem 1 &\moem&\moem\\  0 &\moem 1 &\moem 1 &\moem 1 &\moem 0     \end{array} } \rri   }} $ \\
&&$2$&$3$&${     {  \! S_1'\!=\!\lle {\arw \begin{array}{ccccc}&\moem&\moem 1 &\moem&\moem \\ &\moem&\moem 1 &\moem&\moem\\  1 &\moem 1 &\moem 1 &\moem 0 &\moem 0     \end{array} } \rri   }}$&$ \!{    S_2'\!=\!\lle {\arw \begin{array}{ccccc}&\moem&\moem 0 &\moem&\moem \\ &\moem&\moem 0 &\moem&\moem\\  0 &\moem 0 &\moem 1 &\moem 1 &\moem 0     \end{array} } \rri   }$&$ {     {  \! S_3'\!=\!\lle {\arw \begin{array}{ccccc}&\moem&\moem 0 &\moem&\moem \\ &\moem&\moem 1 &\moem&\moem\\ 0 &\moem 1 &\moem 1 &\moem 1 &\moem 1      \end{array} } \rri   }}$ \\
&&$3$&$2$&${     {  \! S_1''\!=\!\lle {\arw \begin{array}{ccccc}&\moem&\moem 0 &\moem&\moem \\ &\moem&\moem 1 &\moem&\moem\\   1 &\moem 1 &\moem 1 &\moem 1 &\moem 0    \end{array} } \rri   }}$&$ {     {  \! S_2''\!=\!\lle {\arw \begin{array}{ccccc}&\moem&\moem 1 &\moem&\moem \\ &\moem&\moem 1 &\moem&\moem\\  0 &\moem 1 &\moem 2 &\moem 1 &\moem 1     \end{array} } \rri   }}$& \\\hline
&&&&&&\vspace*{-0.615cm}\\
&$\!\!P(3)\!\!$&$1$&$3$&${     {  \! S_1\!=\!\lle {\arw \begin{array}{ccccc}&\moem&\moem 0 &\moem&\moem \\ &\moem&\moem 0 &\moem&\moem\\  1 &\moem 1 &\moem 1 &\moem 1 &\moem 0     \end{array} } \rri   }}$&$ {     {  \! S_2\!=\!\lle {\arw \begin{array}{ccccc}&\moem&\moem 0 &\moem&\moem \\ &\moem&\moem 1 &\moem&\moem\\  0 &\moem 0 &\moem 1 &\moem 1 &\moem 1     \end{array} } \rri   }}$&$ {     {  \! S_3\!=\!\lle {\arw \begin{array}{ccccc}&\moem&\moem 1 &\moem&\moem \\ &\moem&\moem 1 &\moem&\moem\\   0 &\moem 1 &\moem 1 &\moem 0 &\moem 0    \end{array} } \rri   }}$ \\
&&$2$&$3$&${     {  \! S_1'\!=\!\lle {\arw \begin{array}{ccccc}&\moem&\moem 0 &\moem&\moem \\ &\moem&\moem 1 &\moem&\moem\\  1 &\moem 1 &\moem 1 &\moem 0 &\moem 0     \end{array} } \rri   }}$&$ {     {  \! S_2'\!=\!\lle {\arw \begin{array}{ccccc}&\moem&\moem 1 &\moem&\moem \\ &\moem&\moem 1 &\moem&\moem\\ 0 &\moem 0 &\moem 1 &\moem 1 &\moem 0   \end{array} } \rri   }}$&$ {     {  \! S_3'\!=\!\lle {\arw \begin{array}{ccccc} &\moem&\moem 0 &\moem&\moem \\ &\moem&\moem 0 &\moem&\moem\\ 0 &\moem 1 &\moem 1 &\moem 1 &\moem 1     \end{array} } \rri   }}$ \\
&&$3$&$2$&${     {  \! S_1''\!=\!\lle {\arw \begin{array}{ccccc}&\moem&\moem 0 &\moem&\moem \\ &\moem&\moem 1 &\moem&\moem\\ 0 &\moem 1 &\moem 1 &\moem 1 &\moem 0      \end{array} } \rri   }}$&$ {     {  \! S_2''\!=\!\lle {\arw \begin{array}{ccccc}&\moem&\moem 1 &\moem&\moem \\ &\moem&\moem 1 &\moem&\moem\\   1 &\moem 1 &\moem 2 &\moem 1 &\moem 1    \end{array} } \rri   }} $& \\
\hline\hline
%
%
$\!\!\!\tilde{E}_7\!\!\!$&$\!\!P(2)\!\!$&$1$&$4$&$\!{    S_1\!=\!\lle {\arw \begin{array}{ccccccc} &\moem&\moem& \moem 0    &\moem&\moem&\moem\\0&\moem 0&\moem 0&\moem 1&\moem 1&\moem 0&\moem 0\end{array}    }\rri}$&$ {     {  \!S_2\!=\!\lle {\arw \begin{array}{ccccccc} &\moem&\moem  &\moem 1   &\moem&\moem&\moem\\0&\moem 1&\moem 1&\moem 1&\moem 1&\moem 1&\moem 0\end{array}    }\rri}}$&$ {     {  \!S_3\!=\!\lle {\arw \begin{array}{ccccccc} &\moem&\moem  &\moem 0   &\moem&\moem&\moem\\1&\moem 1&\moem 1&\moem 1&\moem 1&\moem 1&\moem 1\end{array}    }\rri}}$\\
&&&&&&$\!{    S_4\!=\!\lle {\arw \begin{array}{ccccccc} &\moem&\moem  &\moem 1   &\moem&\moem&\moem\\0&\moem 0&\moem 1&\moem 1&\moem 0&\moem 0&\moem 0\end{array}  }\rri} $\\
&&$2$&$3$&$\!{    S_1'\!=\!\lle {\arw \begin{array}{ccccccc} &\moem&\moem   &\moem 0  &\moem&\moem&\moem\\0&\moem 0&\moem 1&\moem 1&\moem 1&\moem 1&\moem 0\end{array} }\rri}$&$ {     {  \!S_2'\!=\!\lle {\arw \begin{array}{ccccccc} &\moem&\moem  &\moem 1   &\moem&\moem&\moem\\0&\moem 1&\moem 1&\moem 2&\moem 1&\moem 1&\moem 1\end{array}    }\rri}}$&$ {     {  \!S_3'\!=\!\lle {\arw \begin{array}{ccccccc} &\moem&\moem &\moem 1    &\moem&\moem&\moem\\1&\moem 1&\moem 1&\moem 1&\moem 1&\moem 0&\moem 0\end{array}    }\rri}}$\\
&&$3$&$2$&$ {     {  \!S_1''\!=\!\lle {\arw \begin{array}{ccccccc} &\moem&\moem  &\moem 1   &\moem&\moem&\moem\\0&\moem 1&\moem 2&\moem 2&\moem 2&\moem 1&\moem 1\end{array}    }\rri}}$&$ {     {  \!S_2''\!=\!\lle {\arw \begin{array}{ccccccc} &\moem&\moem  &\moem 1   &\moem&\moem&\moem\\1&\moem 1&\moem 1&\moem 2&\moem 1&\moem 1&\moem 0\end{array}    }\rri}}$&\\\hline
&$\!\!P(3)\!\!$&$1$&$4$&$\!{    S_1\!=\!\lle {\arw \begin{array}{ccccccc} &\moem&\moem  &\moem 0   &\moem&\moem&\moem\\0&\moem 0&\moem 0&\moem 1&\moem 1&\moem 0&\moem 0\end{array}  }\rri }$&$ {     {  \! S_2\!=\!\lle {\arw \begin{array}{ccccccc} &\moem&\moem   &\moem 1  &\moem&\moem&\moem\\0&\moem 0&\moem 1&\moem 1&\moem 1&\moem 1&\moem 0\end{array}    }\rri}}$&$ {     {  \! S_3\!=\!\lle {\arw \begin{array}{ccccccc} &\moem&\moem  &\moem 0   &\moem&\moem&\moem\\0&\moem 1&\moem 1&\moem 1&\moem 1&\moem 1&\moem 1\end{array}    }\rri}}$\\
&&&&&&$ {     {  \! S_4\!=\!\lle {\arw \begin{array}{ccccccc} &\moem&\moem  &\moem 1   &\moem&\moem&\moem\\1&\moem 1&\moem 1&\moem 1&\moem 0&\moem 0&\moem 0\end{array}    }\rri}}$ \\
&&$2$&$3$&$ {     {  \! S_1'\!=\!\lle {\arw \begin{array}{ccccccc} &\moem&\moem  &\moem 0   &\moem&\moem&\moem\\1&\moem 1&\moem 1&\moem 1&\moem 1&\moem 1&\moem 0\end{array}    }\rri}}$&$ {     {  \! S_2'\!=\!\lle {\arw \begin{array}{ccccccc} &\moem&\moem   &\moem 1  &\moem&\moem&\moem\\0&\moem 0&\moem 1&\moem 2&\moem 1&\moem 1&\moem 1\end{array}    }\rri}}$&$ {     {  \!S_3'\!=\!\lle {\arw \begin{array}{ccccccc} &\moem&\moem    &\moem 1 &\moem&\moem&\moem\\0&\moem 1&\moem 1&\moem 1&\moem 1&\moem 0&\moem 0\end{array}    }\rri}}$ \\
&&$3$&$2$&$ {     {  \! S_1''\!=\!\lle {\arw \begin{array}{ccccccc} &\moem&\moem  &\moem 1   &\moem&\moem&\moem\\0&\moem 1&\moem 1&\moem 2&\moem 1&\moem 1&\moem 0\end{array}    }\rri}}$&$ {     {  \! S_2''\!=\!\lle {\arw \begin{array}{ccccccc} &\moem&\moem   &\moem 1  &\moem&\moem&\moem\\1&\moem 1&\moem 2&\moem 2&\moem 2&\moem 1&\moem 1\end{array}    }\rri}}$& \\\hline
&$\!\!P(4)\!\!$&$1$&$4$&$ {     {  \! S_1\!=\!\lle {\arw \begin{array}{ccccccc} &\moem&\moem  &\moem 1   &\moem&\moem&\moem\\0&\moem 0&\moem 0&\moem 1&\moem 1&\moem 1&\moem 0\end{array}    }\rri}}$&$ {     {  \!S_2\!=\!\lle {\arw \begin{array}{ccccccc} &\moem&\moem  &\moem 0   &\moem&\moem&\moem\\0&\moem 0&\moem 1&\moem 1&\moem 1&\moem 1&\moem 1\end{array}    }\rri}}$&$ {     {  \!S_3\!=\!\lle {\arw \begin{array}{ccccccc} &\moem&\moem    &\moem 1 &\moem&\moem&\moem\\0&\moem 1&\moem 1&\moem 1&\moem 0&\moem 0&\moem 0\end{array}    }\rri}}$\\
&&&&&&$ {     {  \!S_4\!=\!\lle {\arw \begin{array}{ccccccc} &\moem&\moem  &\moem 0   &\moem&\moem&\moem\\1&\moem 1&\moem 1&\moem 1&\moem 1&\moem 0&\moem 0\end{array}    }\rri}} $\\
&&$2$&$3$&$ {     {  \! S_1'\!=\!\lle {\arw \begin{array}{ccccccc} &\moem&\moem  &\moem 1   &\moem&\moem&\moem\\0&\moem 0&\moem 1&\moem 1&\moem 1&\moem 0&\moem 0\end{array}    }\rri}}$&$ {     {  \!S_2'\!=\!\lle {\arw \begin{array}{ccccccc} &\moem&\moem  &\moem 0   &\moem&\moem&\moem\\0&\moem 1&\moem 1&\moem 1&\moem 1&\moem 1&\moem 0\end{array}    }\rri}}$&$ {     {  \! S_3'\!=\!\lle {\arw \begin{array}{ccccccc} &\moem&\moem &\moem 1    &\moem&\moem&\moem\\1&\moem 1&\moem 1&\moem 2&\moem 1&\moem 1&\moem 1\end{array}    }\rri}}$ \\
&&$3$&$2$&$ {     {  \! S_1''\!=\!\lle {\arw \begin{array}{ccccccc} &\moem&\moem   &\moem 1  &\moem&\moem&\moem\\0&\moem 1&\moem 1&\moem 2&\moem 2&\moem 1&\moem 1\end{array}    }\rri}}$&$ {     {  \! S_2''\!=\!\lle {\arw \begin{array}{ccccccc} &\moem&\moem  &\moem 1   &\moem&\moem&\moem\\1&\moem 1&\moem 2&\moem 2&\moem 1&\moem 1&\moem 0\end{array}    }\rri}}$& \\\hline 
&$\!\!P(8)\!\!$&$1$&$4$&$\!{    S_1\!=\!\lle {\arw \begin{array}{ccccccc} &\moem&\moem   &\moem 0  &\moem&\moem&\moem\\0&\moem 0&\moem 0&\moem 1&\moem 1&\moem 1&\moem 0\end{array} } \rri }$&$ {     {  \! S_2\!=\!\lle {\arw \begin{array}{ccccccc} &\moem&\moem &\moem 1    &\moem&\moem&\moem\\0&\moem 0&\moem 1&\moem 1&\moem 1&\moem 1&\moem 1\end{array}    }\rri}}$&$\!{    S_3\!=\!\lle {\arw \begin{array}{ccccccc} &\moem&\moem   &\moem 0  &\moem&\moem&\moem\\0&\moem 1&\moem 1&\moem 1&\moem 0&\moem 0&\moem 0\end{array}  }\rri}$\\
&&&&&&$ {     {  \! S_4\!=\!\lle {\arw \begin{array}{ccccccc} &\moem&\moem  &\moem 1   &\moem&\moem&\moem\\1&\moem 1&\moem 1&\moem 1&\moem 1&\moem 0&\moem 0\end{array}    }\rri}} $\\
&&$2$&$3$&$\!{    S_1'\!=\!\lle {\arw \begin{array}{ccccccc} &\moem&\moem  &\moem 0   &\moem&\moem&\moem\\0&\moem 0&\moem 1&\moem 1&\moem 1&\moem 0&\moem 0\end{array}  }\rri}$&$ {     {  \! S_2'\!=\!\lle {\arw \begin{array}{ccccccc} &\moem&\moem  &\moem 1   &\moem&\moem&\moem\\0&\moem 1&\moem 1&\moem 1&\moem 1&\moem 1&\moem 0\end{array}    }\rri}}$&$ {     {  \! S_3'\!=\!\lle {\arw \begin{array}{ccccccc} &\moem&\moem  &\moem 1   &\moem&\moem&\moem\\1&\moem 1&\moem 1&\moem 2&\moem 1&\moem 1&\moem 1\end{array}    }\rri}}$ \\
&&$3$&$2$&$ {     {  \! S_1''\!=\!\lle {\arw \begin{array}{ccccccc} &\moem&\moem  &\moem 1   &\moem&\moem&\moem\\0&\moem 1&\moem 1&\moem 2&\moem 2&\moem 1&\moem 1\end{array}    }\rri}}$&$ {     {  \! S_2''\!=\!\lle {\arw \begin{array}{ccccccc} &\moem&\moem  &\moem 1   &\moem&\moem&\moem\\1&\moem 1&\moem 2&\moem 2&\moem 1&\moem 1&\moem 0\end{array}    }\rri}}$& \\\hline\hline
%
%
$\!\!\!\tilde{E}_8\!\!\!$&$\!\!P(1)\!\!$&$1$&$5$&$\!{    \! S_1\!=\! \lle {\arw \begin{array}{cccccccc} &\moem&\moem 0     &\moem&\moem&\moem  &\moem&\moem \\      0 &\moem 0 &\moem 1 &\moem 1 &\moem 1 &\moem 0 &\moem 0 &\moem 0 \end{array}    }\rri   }$&$ {     {  \!S_2\!=\! \lle  {\arw \begin{array}{cccccccc} &\moem&\moem 1      &\moem&\moem&\moem  &\moem&\moem \\ 1&\moem 1 &\moem 1 &\moem 1 &\moem 1 &\moem 1 &\moem 0 &\moem 0 \end{array}    }\rri   }} $&$ \!{    \! S_3\!=\! \lle {\arw \begin{array}{cccccccc} &\moem&\moem 0    &\moem&\moem&\moem  &\moem&\moem \\       0 &\moem 1 &\moem 1 &\moem 1 &\moem 1 &\moem 1 &\moem 1 &\moem 0 \end{array}    }\rri   }$\\
&&&&&$ {     {  \! S_4\!=\!\lle    {\arw \begin{array}{cccccccc} &\moem&\moem 1       &\moem&\moem&\moem  &\moem&\moem \\       1 &\moem 1 &\moem 2 &\moem 1 &\moem 1 &\moem 1 &\moem 1 &\moem 1 \end{array}    }\rri   }}$&$ \!{    \!  S_5\!=\! \lle {\arw \begin{array}{cccccccc} &\moem&\moem 1     &\moem&\moem&\moem  &\moem&\moem \\       0 &\moem 1 &\moem 1 &\moem 1 &\moem 0 &\moem 0 &\moem 0 &\moem 0 \end{array}    }\rri   }  $\\
&&$2$&$3$&$\!{    \!  S_1'\!=\!\lle {\arw \begin{array}{cccccccc} &\moem&\moem 1     &\moem&\moem&\moem  &\moem&\moem \\       0 &\moem 1 &\moem 2 &\moem 1 &\moem 1 &\moem 1 &\moem 0 &\moem 0 \end{array}    }\rri   }$&$ {     {  \!S_2'\!=\!\lle    {\arw \begin{array}{cccccccc} &\moem&\moem 1       &\moem&\moem&\moem  &\moem&\moem \\       1 &\moem 1 &\moem 2 &\moem 2 &\moem 1 &\moem 1 &\moem 1 &\moem 0 \end{array}    }\rri   }}$&${     {  \! S_3'\!=\!\lle    {\arw \begin{array}{cccccccc} &\moem&\moem 1     &\moem&\moem&\moem  &\moem&\moem \\       1 &\moem 2 &\moem 2 &\moem 2 &\moem 2 &\moem 1 &\moem 1 &\moem 1 \end{array}    }\rri   }}  $\\
&&$3$&$2$&${     {  \!S_1''\!=\!\lle    {\arw \begin{array}{cccccccc} &\moem&\moem 1     &\moem&\moem&\moem  &\moem&\moem \\       1 &\moem 2 &\moem 3 &\moem 3 &\moem 2 &\moem 2 &\moem 1 &\moem 1 \end{array}    }\rri   }}$&$ {     {  \!S_2''\!=\!\lle    {\arw \begin{array}{cccccccc} &\moem&\moem 2       &\moem&\moem&\moem  &\moem&\moem \\       1 &\moem 2 &\moem 3 &\moem 2 &\moem 2 &\moem 1 &\moem 1 &\moem 0 \end{array}    }\rri   }}$&$  $\\\hline
&$\!\!P(2)\!\!$&$1$&$5$&$\!{    \!  S_1\!=\!\lle {\arw \begin{array}{cccccccc} &\moem&\moem 0    &\moem&\moem&\moem  &\moem&\moem \\       0 &\moem 0 &\moem 1 &\moem 1 &\moem 1 &\moem 0 &\moem 0 &\moem 0 \end{array}    }\rri   }$&$ {     {  \! S_2\!=\!\lle    {\arw \begin{array}{cccccccc} &\moem&\moem1      &\moem&\moem&\moem  &\moem&\moem \\       0 &\moem 1 &\moem 1 &\moem 1 &\moem 1 &\moem 1 &\moem 0 &\moem 0 \end{array}    }\rri   }}$&$ {     {  \!S_3\!=\!\lle    {\arw \begin{array}{cccccccc} &\moem&\moem 0     &\moem&\moem&\moem  &\moem&\moem \\       1 &\moem 1 &\moem 1 &\moem 1 &\moem 1 &\moem 1 &\moem 1 &\moem 0 \end{array}    }\rri   }}$\\
&&&&&${     {  \!S_4\!=\!\lle    {\arw \begin{array}{cccccccc} &\moem&\moem 1     &\moem&\moem&\moem  &\moem&\moem \\       0 &\moem 1 &\moem 2 &\moem 1 &\moem 1 &\moem 1 &\moem 1 &\moem 1 \end{array}    }\rri   }}$&$ 
{     {  \! S_5\!=\!\lle    {\arw \begin{array}{cccccccc} &\moem&\moem 1     &\moem&\moem&\moem  &\moem&\moem \\       1 &\moem 1 &\moem 1 &\moem 1 &\moem 0 &\moem 0 &\moem 0 &\moem 0 \end{array}    }\rri   }}  $\\
&&$2$&$3$&${     {  \!S_1'\!=\!\lle    {\arw \begin{array}{cccccccc} &\moem&\moem 1     &\moem&\moem&\moem  &\moem&\moem \\       1 &\moem 2 &\moem 2 &\moem 2 &\moem 2 &\moem 1 &\moem 1 &\moem 1 \end{array}    }\rri   }}$&$ {     {  \!S_2'\!=\!\lle    {\arw \begin{array}{cccccccc} &\moem&\moem 1      &\moem&\moem&\moem  &\moem&\moem \\       1 &\moem 1 &\moem 2 &\moem 1 &\moem 1 &\moem 1 &\moem 0 &\moem 0 \end{array}    }\rri   }}$&${     {  \!S_3'\!=\!\lle    {\arw \begin{array}{cccccccc} &\moem&\moem 1     &\moem&\moem&\moem  &\moem&\moem \\       0 &\moem 1 &\moem 2 &\moem 2 &\moem 1 &\moem 1 &\moem 1 &\moem 0 \end{array}    }\rri   }}  $\\
&&$3$&$2$&${     {  \!S_1''\!=\!\lle    {\arw \begin{array}{cccccccc} &\moem&\moem 1     &\moem&\moem&\moem  &\moem&\moem \\       1 &\moem 2 &\moem 3 &\moem 3 &\moem 2 &\moem 2 &\moem 1 &\moem 1 \end{array}    }\rri   }}$&$ {     {  \!S_2''\!=\!\lle    {\arw \begin{array}{cccccccc} &\moem&\moem 2      &\moem&\moem&\moem  &\moem&\moem \\       1 &\moem 2 &\moem 3 &\moem 2 &\moem 2 &\moem 1 &\moem 1 &\moem 0 \end{array}    }\rri   }}$&$  $\\\hline
&$\!\!P(3)\!\!$&$1$&$5$&${     {  \!S_1\!=\!\lle    {\arw \begin{array}{cccccccc} &\moem&\moem 1     &\moem&\moem&\moem  &\moem&\moem \\       0 &\moem 0 &\moem 1 &\moem 1 &\moem 1 &\moem 1 &\moem 0 &\moem 0 \end{array}    }\rri   }}$&$ {     {  \!S_2\!=\!\lle    {\arw \begin{array}{cccccccc} &\moem&\moem 0    &\moem&\moem&\moem  &\moem&\moem \\       0 &\moem 1 &\moem 1 &\moem 1 &\moem 1 &\moem 1 &\moem 1 &\moem 0 \end{array}    }\rri   }}$&$ {     {  \!S_3\!=\!\lle    {\arw \begin{array}{cccccccc} &\moem&\moem 1    &\moem&\moem&\moem  &\moem&\moem \\       1 &\moem 1 &\moem 2 &\moem 1 &\moem 1 &\moem 1 &\moem 1 &\moem 1 \end{array}    }\rri   }}$\\
&&&&&$ {     {  \!S_4\!=\!\lle    {\arw \begin{array}{cccccccc} &\moem&\moem 1       &\moem&\moem&\moem  &\moem&\moem \\       0 &\moem 1 &\moem 1 &\moem 1 &\moem 0 &\moem 0 &\moem 0 &\moem 0 \end{array}    }\rri   }}$&${     {  \!S_5\!=\!\lle    {\arw \begin{array}{cccccccc} &\moem&\moem 0    &\moem&\moem&\moem  &\moem&\moem \\       1 &\moem 1 &\moem 1 &\moem 1 &\moem 1 &\moem 0 &\moem 0 &\moem 0 \end{array}    }\rri   }} $\\
&&$2$&$3$&${     {  \!S_1'\!=\!\lle    {\arw \begin{array}{cccccccc} &\moem&\moem 1    &\moem&\moem&\moem  &\moem&\moem \\       1 &\moem 2 &\moem 2 &\moem 1 &\moem 1 &\moem 1 &\moem 0 &\moem 0 \end{array}    }\rri   }}$&$ {     {  \!S_2'\!=\!\lle    {\arw \begin{array}{cccccccc} &\moem&\moem 1    &\moem&\moem&\moem  &\moem&\moem \\       1 &\moem 1 &\moem 2 &\moem 2 &\moem 1 &\moem 1 &\moem 1 &\moem 0 \end{array}    }\rri   }}$&${     {  \!S_3'\!=\!\lle    {\arw \begin{array}{cccccccc} &\moem&\moem 1     &\moem&\moem&\moem  &\moem&\moem \\       0 &\moem 1 &\moem 2 &\moem 2 &\moem 2 &\moem 1 &\moem 1 &\moem 1 \end{array}    }\rri   }}  $\\
&&$3$&$2$&${     {  \!S_1''\!=\!\lle    {\arw \begin{array}{cccccccc} &\moem&\moem 1    &\moem&\moem&\moem  &\moem&\moem \\       1 &\moem 2 &\moem 3 &\moem 3 &\moem 2 &\moem 2 &\moem 1 &\moem 1 \end{array}    }\rri   }}$&$ {     {  \!S_2''\!=\!\lle    {\arw \begin{array}{cccccccc} &\moem&\moem 2     &\moem&\moem&\moem  &\moem&\moem \\       1 &\moem 2 &\moem 3 &\moem 2 &\moem 2 &\moem 1 &\moem 1 &\moem 0 \end{array}    }\rri   }}$&$  $\\\hline
&$\!\!P(4)\!\!$&$1$&$5$&${     {  \!S_1\!=\!\lle    {\arw \begin{array}{cccccccc} &\moem&\moem 1    &\moem&\moem&\moem  &\moem&\moem \\       0 &\moem 0 &\moem 1 &\moem 1 &\moem 1 &\moem 1 &\moem 1 &\moem 0 \end{array}    }\rri   }}$&$ {     {  \!S_2\!=\!\lle    {\arw \begin{array}{cccccccc} &\moem&\moem 0      &\moem&\moem&\moem  &\moem&\moem \\       0 &\moem 1 &\moem 1 &\moem 1 &\moem 1 &\moem 1 &\moem 1 &\moem 1 \end{array}    }\rri   }}$&$ {     {  \!S_3\!=\!\lle    {\arw \begin{array}{cccccccc} &\moem&\moem 1     &\moem&\moem&\moem  &\moem&\moem \\       1 &\moem 1 &\moem 2 &\moem 1 &\moem 0 &\moem 0 &\moem 0 &\moem 0 \end{array}    }\rri   }}$\\
&&&&& ${     {  \!S_4\!=\!\lle    {\arw \begin{array}{cccccccc} &\moem&\moem 1      &\moem&\moem&\moem  &\moem&\moem \\       0 &\moem 1 &\moem 1 &\moem 1 &\moem 1 &\moem 0 &\moem 0 &\moem 0 \end{array}    }\rri   }}$&${     {  \!S_5\!=\!\lle    {\arw \begin{array}{cccccccc} &\moem&\moem 0     &\moem&\moem&\moem  &\moem&\moem \\       1 &\moem 1 &\moem 1 &\moem 1 &\moem 1 &\moem 1 &\moem 0 &\moem 0 \end{array}    }\rri   }}  $\\
&&$2$&$3$&${     {  \!S_1'\!=\!\lle    {\arw \begin{array}{cccccccc} &\moem&\moem 1      &\moem&\moem&\moem  &\moem&\moem \\       1 &\moem 1 &\moem 2 &\moem 2 &\moem 2 &\moem 1 &\moem 1 &\moem 1 \end{array}    }\rri   }}$&$ {     {  \!S_2'\!=\!\lle    {\arw \begin{array}{cccccccc} &\moem&\moem 1      &\moem&\moem&\moem  &\moem&\moem \\       0 &\moem 1 &\moem 2 &\moem 1 &\moem 1 &\moem 1 &\moem 0 &\moem 0 \end{array}    }\rri   }}$&${     {  \!S_3'\!=\!\lle    {\arw \begin{array}{cccccccc} &\moem&\moem 1     &\moem&\moem&\moem  &\moem&\moem \\       1 &\moem 2 &\moem 2 &\moem 2 &\moem 1 &\moem 1 &\moem 1 &\moem 0 \end{array}    }\rri   }}  $\\
&&$3$&$2$&${     {  \!S_1''\!=\!\lle    {\arw \begin{array}{cccccccc} &\moem&\moem 2     &\moem&\moem&\moem  &\moem&\moem \\       1 &\moem 2 &\moem 3 &\moem 3 &\moem 2 &\moem 2 &\moem 1 &\moem 1 \end{array}    }\rri   }}$&$ {     {  \!S_2''\!=\!\lle    {\arw \begin{array}{cccccccc} &\moem&\moem 1      &\moem&\moem&\moem  &\moem&\moem \\       1 &\moem 2 &\moem 3 &\moem 2 &\moem 2 &\moem 1 &\moem 1 &\moem 0 \end{array}    }\rri   }}$&$  $\\\hline
&$\!\!P(5)\!\!$&$1$&$5$&${     {  \!S_1\!=\!\lle    {\arw \begin{array}{cccccccc} &\moem&\moem 1     &\moem&\moem&\moem  &\moem&\moem \\       0 &\moem 0 &\moem 1 &\moem 1 &\moem 1 &\moem 1 &\moem 1 &\moem 1 \end{array}    }\rri   }}$&$ \!{    \! S_2\!=\! \lle {\arw \begin{array}{cccccccc} &\moem&\moem 0      &\moem&\moem&\moem  &\moem&\moem \\       0 &\moem 1 &\moem 1 &\moem 1 &\moem 0 &\moem 0 &\moem 0 &\moem 0 \end{array}    }\rri   }$&${     {  \!S_3\!=\!\lle    {\arw \begin{array}{cccccccc} &\moem&\moem 1    &\moem&\moem&\moem  &\moem&\moem \\       1 &\moem 1 &\moem 2 &\moem 1 &\moem 1 &\moem 0 &\moem 0 &\moem 0 \end{array}    }\rri   }}$\\
&&&&& ${     {  \! S_4\!=\!\lle    {\arw \begin{array}{cccccccc} &\moem&\moem 1     &\moem&\moem&\moem  &\moem&\moem \\       0 &\moem 1 &\moem 1 &\moem 1 &\moem 1 &\moem 1 &\moem 0 &\moem 0 \end{array}    }\rri   }}$&${     {  \!S_5\!=\!\lle    {\arw \begin{array}{cccccccc} &\moem&\moem 0    &\moem&\moem&\moem  &\moem&\moem \\       1 &\moem 1 &\moem 1 &\moem 1 &\moem 1 &\moem 1 &\moem 1 &\moem 0 \end{array}    }\rri   }}  $\\
&&$2$&$3$&${     {  \!S_1'\!=\!\lle    {\arw \begin{array}{cccccccc} &\moem&\moem 1      &\moem&\moem&\moem  &\moem&\moem \\       1 &\moem 1 &\moem 2 &\moem 2 &\moem 1 &\moem 1 &\moem 0 &\moem 0 \end{array}    }\rri   }}$&$ {     {  \!S_2'\!=\!\lle    {\arw \begin{array}{cccccccc} &\moem&\moem 1      &\moem&\moem&\moem  &\moem&\moem \\       0 &\moem 1 &\moem 2 &\moem 1 &\moem 1 &\moem 1 &\moem 1 &\moem 0 \end{array}    }\rri   }}$&$ {     {  \!S_3'\!=\!\lle    {\arw \begin{array}{cccccccc} &\moem&\moem 1    &\moem&\moem&\moem  &\moem&\moem \\       1 &\moem 2 &\moem 2 &\moem 2 &\moem 2 &\moem 1 &\moem 1 &\moem 1 \end{array}    }\rri   }}  $\\
&&$3$&$2$&${     {  \!S_1''\!=\!\lle    {\arw \begin{array}{cccccccc} &\moem&\moem 2     &\moem&\moem&\moem  &\moem&\moem \\       1 &\moem 2 &\moem 3 &\moem 3 &\moem 2 &\moem 1 &\moem 1 &\moem 0 \end{array}    }\rri   }}$&$ {     {  \!S_2''\!=\!\lle    {\arw \begin{array}{cccccccc} &\moem&\moem 1      &\moem&\moem&\moem  &\moem&\moem \\       1 &\moem 2 &\moem 3 &\moem 2 &\moem 2 &\moem 2 &\moem 1 &\moem 1 \end{array}    }\rri   }}$&$  $\\\hline
&$\!\!P(6)\!\!$&$1$&$5$&$\!{    \! S_1\!=\! \lle {\arw \begin{array}{cccccccc} &\moem&\moem 0     &\moem&\moem&\moem  &\moem&\moem \\       0 &\moem 1 &\moem 1 &\moem 1 &\moem 0 &\moem 0 &\moem 0 &\moem 0 \end{array}    }\rri   }$&$ {     {  \!S_2\!=\!\lle    {\arw \begin{array}{cccccccc} &\moem&\moem 1      &\moem&\moem&\moem  &\moem&\moem \\       1 &\moem 1 &\moem 2 &\moem 1 &\moem 1 &\moem 1 &\moem 0 &\moem 0 \end{array}    }\rri   }}$&$ {     {  \!S_3\!=\!\lle    {\arw \begin{array}{cccccccc} &\moem&\moem 1      &\moem&\moem&\moem  &\moem&\moem \\       0 &\moem 1 &\moem 1 &\moem 1 &\moem 1 &\moem 1 &\moem 1 &\moem 0\end{array}    }\rri   }}$\\
&&&&& ${     {  \!S_4\!=\!\lle    {\arw \begin{array}{cccccccc} &\moem&\moem 0     &\moem&\moem&\moem  &\moem&\moem \\       1 &\moem 1 &\moem 1 &\moem 1 &\moem 1 &\moem 1 &\moem 1 &\moem 1 \end{array}    }\rri   }}$&$\!{    \!  S_5\!=\!\lle {\arw \begin{array}{cccccccc} &\moem&\moem 1    &\moem&\moem&\moem  &\moem&\moem \\       0 &\moem 0 &\moem 1 &\moem 1 &\moem 1 &\moem 0 &\moem 0 &\moem 0 \end{array}    }\rri   }  $\\
&&$2$&$3$&${     {  \!S_1'\!=\!\lle    {\arw \begin{array}{cccccccc} &\moem&\moem 1     &\moem&\moem&\moem  &\moem&\moem \\       1 &\moem 1 &\moem 2 &\moem 2 &\moem 1 &\moem 1 &\moem 1 &\moem 0 \end{array}    }\rri   }}$&$ {     {  \!S_2'\!=\!\lle    {\arw \begin{array}{cccccccc} &\moem&\moem 1      &\moem&\moem&\moem  &\moem&\moem \\       0 &\moem 1 &\moem 2 &\moem 1 &\moem 1 &\moem 1 &\moem 1 &\moem 1 \end{array}    }\rri   }}$&${     {  \!S_3'\!=\!\lle    {\arw \begin{array}{cccccccc} &\moem&\moem 1    &\moem&\moem&\moem  &\moem&\moem \\       1 &\moem 2 &\moem 2 &\moem 2 &\moem 2 &\moem 1 &\moem 0 &\moem 0 \end{array}    }\rri   }} $\\
&&$3$&$2$&${     {  \!S_1''\!=\!\lle    {\arw \begin{array}{cccccccc} &\moem&\moem 2     &\moem&\moem&\moem  &\moem&\moem \\       1 &\moem 2 &\moem 3 &\moem 3 &\moem 2 &\moem 2 &\moem 1 &\moem 1 \end{array}    }\rri   }}$&$ {     {  \!S_2''\!=\!\lle    {\arw \begin{array}{cccccccc} &\moem&\moem 1      &\moem&\moem&\moem  &\moem&\moem \\       1 &\moem 2 &\moem 3 &\moem 2 &\moem 2 &\moem 1 &\moem 1 &\moem 0 \end{array}    }\rri   }}$&  \\\hline
&$\!\!P(7)\!\!$&$1$&$5$&$\!{    \!  S_1\!=\!\lle {\arw \begin{array}{cccccccc} &\moem&\moem 0     &\moem&\moem&\moem  &\moem&\moem \\       0 &\moem 1 &\moem 1 &\moem 1 &\moem 0 &\moem 0 &\moem 0 &\moem 0 \end{array}    }\rri   }$&$ {     {  \!S_2\!=\!\lle    {\arw \begin{array}{cccccccc} &\moem&\moem 1     &\moem&\moem&\moem  &\moem&\moem \\       1 &\moem 1 &\moem 2 &\moem 1 &\moem 1 &\moem 1 &\moem 1 &\moem 0 \end{array}    }\rri   }}$&${     {  \!S_3\!=\!\lle    {\arw \begin{array}{cccccccc} &\moem&\moem 1      &\moem&\moem&\moem  &\moem&\moem \\       0 &\moem 1 &\moem 1 &\moem 1 &\moem 1 &\moem 1 &\moem 1 &\moem 1 \end{array}    }\rri   }}$\\
&&&&&$ \!{    \!  S_4\!=\!\lle {\arw \begin{array}{cccccccc} &\moem&\moem 0      &\moem&\moem&\moem  &\moem&\moem \\       1 &\moem 1 &\moem 1 &\moem 1 &\moem 1 &\moem 1 &\moem 0 &\moem 0 \end{array}    }\rri   }$&$ \!{    \!  S_5\!=\!\lle {\arw \begin{array}{cccccccc} &\moem&\moem 1    &\moem&\moem&\moem  &\moem&\moem \\       0 &\moem 0 &\moem 1 &\moem 1 &\moem 1 &\moem 0 &\moem 0 &\moem 0 \end{array}    }\rri   }  $\\
&&$2$&$3$&${     {  \!S_1'\!=\!\lle    {\arw \begin{array}{cccccccc} &\moem&\moem 1      &\moem&\moem&\moem  &\moem&\moem \\       1 &\moem 2 &\moem 2 &\moem 2 &\moem 2 &\moem 1 &\moem 1 &\moem 0 \end{array}    }\rri   }}$&$ {     {  \!S_2'\!=\!\lle    {\arw \begin{array}{cccccccc} &\moem&\moem 1       &\moem&\moem&\moem  &\moem&\moem \\       1 &\moem 1 &\moem 2 &\moem 2 &\moem 1 &\moem 1 &\moem 1 &\moem 1 \end{array}    }\rri   }}$&$\!{    \!  S_3'\!=\!\lle {\arw \begin{array}{cccccccc} &\moem&\moem 1      &\moem&\moem&\moem  &\moem&\moem \\       0 &\moem 1 &\moem 2 &\moem 1 &\moem 1 &\moem 1 &\moem 0 &\moem 0 \end{array}    }\rri   }  $\\
&&$3$&$2$&${     {  \!S_1''\!=\!\lle    {\arw \begin{array}{cccccccc} &\moem&\moem 2      &\moem&\moem&\moem  &\moem&\moem \\       1 &\moem 2 &\moem 3 &\moem 3 &\moem 2 &\moem 2 &\moem 1 &\moem 0 \end{array}    }\rri   }}$&$ {     {  \!S_2''\!=\!\lle    {\arw \begin{array}{cccccccc} &\moem&\moem 1      &\moem&\moem&\moem  &\moem&\moem \\       1 &\moem 2 &\moem 3 &\moem 2 &\moem 2 &\moem 1 &\moem 1 &\moem 1 \end{array}    }\rri   }}$&$  $\\\hline
&$\!\!P(9)\!\!$&$1$&$5$&$\!{    \!  S_1\!=\!\lle {\arw \begin{array}{cccccccc} &\moem&\moem 0      &\moem&\moem&\moem  &\moem&\moem \\       0 &\moem 1 &\moem 1 &\moem 1 &\moem 0 &\moem 0 &\moem 0 &\moem 0 \end{array}    }\rri   }$&$ {     {  \!S_2\!=\!\lle    {\arw \begin{array}{cccccccc} &\moem&\moem 1       &\moem&\moem&\moem  &\moem&\moem \\       1 &\moem 1 &\moem 1 &\moem 1 &\moem 1 &\moem 0 &\moem 0 &\moem 0 \end{array}    }\rri   }}$&$ \!{    \! S_3\!=\! \lle {\arw \begin{array}{cccccccc} &\moem&\moem 0     &\moem&\moem&\moem  &\moem&\moem \\       0 &\moem 0 &\moem 1 &\moem 1 &\moem 1 &\moem 1 &\moem 0 &\moem 0 \end{array}    }\rri   }$\\&&&&& ${     {  \!S_4\!=\!\lle    {\arw \begin{array}{cccccccc} &\moem&\moem 1      &\moem&\moem&\moem  &\moem&\moem \\       0 &\moem 1 &\moem 1 &\moem 1 &\moem 1 &\moem 1 &\moem 1 &\moem 0 \end{array}    }\rri   }}$&$ {     {  \!S_5\!=\!\lle    {\arw \begin{array}{cccccccc} &\moem&\moem 1     &\moem&\moem&\moem  &\moem&\moem \\       1 &\moem 1 &\moem 2 &\moem 1 &\moem 1 &\moem 1 &\moem 1 &\moem 1 \end{array}    }\rri   }}  $\\
&&$2$&$3$&${     {  \!S_1'\!=\!\lle    {\arw \begin{array}{cccccccc} &\moem&\moem 1    &\moem&\moem&\moem  &\moem&\moem \\       1 &\moem 1 &\moem 2 &\moem 2 &\moem 1 &\moem 1 &\moem 1 &\moem 0 \end{array}    }\rri   }}$&$ {     {  \!S_2'\!=\!\lle    {\arw \begin{array}{cccccccc} &\moem&\moem 1       &\moem&\moem&\moem  &\moem&\moem \\       0 &\moem 1 &\moem 2 &\moem 2 &\moem 2 &\moem 1 &\moem 1 &\moem 1 \end{array}    }\rri   }}$&$ {     {  \!S_3'\!=\!\lle    {\arw \begin{array}{cccccccc} &\moem&\moem 1      &\moem&\moem&\moem  &\moem&\moem \\       1 &\moem 2 &\moem 2 &\moem 1 &\moem 1 &\moem 1 &\moem 0 &\moem 0 \end{array}    }\rri   }}  $\\
&&$3$&$2$&${     {  \!S_1''\!=\!\lle    {\arw \begin{array}{cccccccc} &\moem&\moem 2      &\moem&\moem&\moem  &\moem&\moem \\       1 &\moem 2 &\moem 3 &\moem 3 &\moem 2 &\moem 2 &\moem 1 &\moem 1 \end{array}    }\rri   }}$&$ {     {  \!S_2''\!=\!\lle    {\arw \begin{array}{cccccccc} &\moem&\moem 1       &\moem&\moem&\moem  &\moem&\moem \\       1 &\moem 2 &\moem 3 &\moem 2 &\moem 2 &\moem 1 &\moem 1 &\moem 0 \end{array}    }\rri   }}$&\\
\end{longtable}

%% file: tables.tex
\renewcommand{\arraystretch}{1.2}\normalsize
\subsection{The minimal deformations from a single tube}\label{mindeglist}
Finally, we give the list of minimal degenerations $M<U\oplus V$ where $U$ is simple projective, $V$ is preinjective and $M=M_\mu$ comes from a non-homogeneous tube. Up to application of the periodicity theoren \ref{periprop} and proposition \ref{indsub}, the list contains all minimal degenerations of this type.
The table includes column by column the following informations:

\vspace{0.25cm}
\begin{tabular}{ll}
(1)&  the type of $Q$;\\
(2) &the only projective simple $U$ (which determines the orientation of $Q$);\\
(3)&  the preinjective indecomposable $V$;\\
(4) &the non-homogeneous tube $\mathcal{T}_{k}$;\\
(5) & the list of the minimal deformations $M_{\,\!k}$ of $U\oplus V$;
\end{tabular}
\begin{longtable}{||c|c|c|c|c||}\hline\hline
$|Q|$&$ U $&$V$&$ k $&$M_{\, k}$ \\\hline
\endhead
\hline\hline\endfoot
$\tilde{D}_8$&$  P(3)  $&$ \tau^{0}I(3) $&$ 1 $&$  S_{3}(1)\oplus S_{4}(1)  $\\
&&$ \tau^{1}I(4) $&$ 1 $&$  S_{2}(2)\oplus S_{4}(1)  $\\
&&$ \tau^{2}I(5) $&$ 1 $&$  S_{1}(3)\oplus S_{4}(1)  $\\
&&$ \tau^{3}I(6) $&$ 1 $&$  S_{4}(1)\oplus S_{6}(4)  $\\
&&$ \tau^{4}I(7) $&$ 1 $&$  S_{4}(1)\oplus S_{5}(5)  $\\
&&$ \tau^{4}I(7) $&$ 2 $&$  S_{1}'(1)\oplus S_{2}'(1)  $\\
&&$ \tau^{4}I(7) $&$ 3 $&$  S_{1}''(1)\oplus S_{2}''(1)  $\\
\hline
&$  P(4)  $&$ \tau^{0}I(4)  $&$ 1 $&$  S_{3}(1)\oplus S_{5}(1)  $\\
&&$ \tau^{1}I(3) $&$ 1 $&$  S_{3}(1)\oplus S_{4}(2)  $\\
&&$ \tau^{1}I(5) $&$ 1 $&$  S_{2}(2)\oplus S_{5}(1)  $\\
&&$ \tau^{2}I(6)  $&$ 1 $&$  S_{1}(3)\oplus S_{5}(1)  $\\
&&$ \tau^{3}I(7)  $&$ 1 $&$  S_{5}(1)\oplus S_{6}(4)  $\\
&&$ \tau^{3}I(6)  $&$ 2 $&$  S_{1}'(1)\oplus S_{2}'(1) $\\
&&$ \tau^{3}I(6)  $&$ 3 $&$  S_{1}''(1)\oplus S_{2}''(1)  $\\

\hline
&$  P(5)  $&$ \tau^{0}I(5)  $&$ 1 $&$  S_{3}(1)\oplus S_{6}(1)  $\\
&&$ \tau^{1}I(4)  $&$ 1 $&$  S_{3}(1)\oplus S_{5}(2)  $\\
&&$ \tau^{2}I(3)  $&$ 1 $&$  S_{3}(1)\oplus S_{4}(3)  $\\
&&$ \tau^{1}I(6)  $&$ 1 $&$  S_{2}(2)\oplus S_{6}(1)  $\\
&&$ \tau^{2}I(7)  $&$ 1 $&$  S_{1}(3)\oplus S_{6}(1)  $\\
&&$ \tau^{2}I(5)  $&$ 2 $&$  S_{1}'(1)\oplus S_{2}'(1)  $\\
&&$ \tau^{2}I(5)  $&$ 3 $&$  S_{1}''(1)\oplus S_{2}''(1)  $\\

\hline\hline
%
%
$\tilde{E}_6$&$  P(2)  $&$ \tau^{2}I(6) $&$ 1 $&$  S_{1}(1)\oplus S_{3}(1)$\\
&&$ \tau^{2}I(2) $&$ 1 $&$  S_{1}(1)\oplus S_{2}(2)$\\
&&$ \tau^{2}I(4) $&$ 2 $&$  S_{1}'(1)\oplus S_{3}'(1) $\\
&&$ \tau^{2}I(2) $&$ 2 $&$  S_{1}'(1)\oplus S_{2}'(2)$\\
&&$ \tau^{2}I(2) $&$ 3 $&$  S_{1}''(1)\oplus S_{2}''(1) $\\
\hline
&$  P(3)  $&$ \tau^{1}I(3) $&$ 1 $&$  S_{1}(1)\oplus S_{2}(1)\oplus S_{3}(1)$\\
&&$ \tau^{1}I(3)$&$ 2 $&$  S_{1}'(1)\oplus S_{2}'(1)\oplus S_{3}'(1) $\\
&&$ \tau^{2}I(3)$&$ 3 $&$  S_{1}''(1)\oplus S_{2}''(1)\oplus S_{2}''(1) $\\
\hline\hline
%
%
$\tilde{E}_7$&$  P(2)  $&$ \tau^{4}I(6) $&$ 1 $&$  S_{2}(1)\oplus S_{3}(1) $\\
&&$ \tau^{5}I(8) $&$ 1 $&$  S_{1}(2)\oplus S_{3}(1)$\\
&&$ \tau^{7}I(6) $&$ 1 $&$  S_{3}(1)\oplus S_{4}(3) $\\
&&$ \tau^{3}I(2) $&$ 2 $&$  S_{2}'(1)\oplus S_{3}'(1)$\\
&&$ \tau^{7}I(6) $&$ 2 $&$  S_{1}'(2)\oplus S_{3}'(1)$\\
&&$ \tau^{7}I(6) $&$ 3 $&$  S_{1}''(1)\oplus S_{2}''(1)$\\
\hline
&$  P(3)  $
&$ \tau^{3}I(5) $&$ 1 $&$  S_{2}(1)\oplus S_{3}(1)\oplus S_{4}(1) $\\
&&$ \tau^{3}I(3) $&$ 1 $&$  S_{1}(2)\oplus S_{3}(1)\oplus S_{4}(1) $\\
&&$ \tau^{5}I(5)$&$ 1 $&$  S_{1}(2)\oplus S_{2}(2)\oplus S_{4}(1) $\\
&&$ \tau^{3}I(3) $&$ 2 $&$  S_{1}'(1)\oplus S_{2}'(1)\oplus S_{3}'(1)$\\
&&$ \tau^{6}I(5)$&$ 3 $&$  S_{1}''(1)\oplus S_{2}''(1)\oplus S_{2}''(1)$\\
\hline
&$  P(4)  $&$ \tau^{2}I(4)$&$ 1 $&$  S_{1}(1)\oplus S_{2}(1)\oplus S_{3}(1)\oplus S_{4}(1)$\\
&&$ \tau^{3}I(4) $&$ 2 $&$  S_{1}'(1)\oplus S_{2}'(1)\oplus S_{3}'(1)\oplus S_{3}'(1) $\\
&&$ \tau^{5}I(4)$&$ 3 $&$  S_{1}''(1)\oplus S_{1}''(1)\oplus S_{2}''(1)\oplus S_{2}''(1) $\\
\hline
&$  P(8)  $
&$ \tau^{2}I(8)$&$ 1 $&$  S_{2}(1)\oplus S_{4}(1)$\\
&&$ \tau^{5}I(2) $&$ 1 $&$  S_{2}(1)\oplus S_{3}(2)$\\
&&$ \tau^{5}I(6) $&$ 1 $&$  S_{1}(2)\oplus S_{4}(1)$\\
&&$ \tau^{3}I(8) $&$ 2 $&$  S_{2}'(1)\oplus S_{3}'(1)$\\
&&$ \tau^{5}I(8) $&$ 2 $&$  S_{1}'(2)\oplus S_{3}'(1)$\\
&&$ \tau^{5}I(8) $&$ 3 $&$  S_{1}''(1)\oplus S_{2}''(1) $\\


\hline\hline
$\tilde{E}_8$&$        P(1)        $
&$ \tau^{5}I(1)  $&$ 1 $&$  S_{2}(1)\oplus S_{4}(1)$\\
&&$ \tau^{8}I(1)  $&$ 1 $&$  S_{1}(2)\oplus S_{4}(1)$\\
&&$ \tau^{14}I(7)  $&$ 1 $&$  S_{4}(1)\oplus S_{5}(3)$\\
&&$ \tau^{11}I(7)  $&$ 1 $&$  S_{2}(1)\oplus S_{3}(2)$\\
&&$ \tau^{9}I(1)  $&$ 2 $&$  S_{2}'(1)\oplus S_{3}'(1)$\\
&&$ \tau^{14}I(1)  $&$ 2 $&$  S_{1}'(2)\oplus S_{3}'(1)$\\
&&$ \tau^{14}I(1)  $&$ 3 $&$  S_{1}''(1)\oplus S_{2}''(1)$\\

\hline
&$        P(2)        $
&$ \tau^{5}I(2)  $&$ 1 $&$  S_{2}(1)\oplus S_{3}(1)\oplus S_{4}(1)\oplus S_{5}(1) $\\
&&$ \tau^{8}I(5)  $&$ 1 $&$  S_{1}(2)\oplus S_{3}(1)\oplus S_{4}(1)\oplus S_{5}(1) $\\
&&$ \tau^{8}I(2)  $&$ 1 $&$  S_{1}(2)\oplus S_{2}(2)\oplus S_{4}(1)\oplus S_{5}(1)$\\
&&$ \tau^{11}I(5) $&$ 1 $&$  S_{1}(2)\oplus S_{2}(2)\oplus S_{3}(2)\oplus S_{5}(1)$\\
&&$ \tau^{9}I(2)  $&$ 2 $&$  S_{1}'(1)\oplus S_{1}'(1)\oplus S_{2}'(1)\oplus S_{3}'(1)$\\
&&$ \tau^{14}I(2)  $&$ 3 $&$  S_{1}''(1)\oplus S_{1}''(1)\oplus S_{2}''(1)\oplus S_{2}''(1) $\\

\hline
&$        P(3)        $&$\tau^{5}I(3) $&$ 1 $&$        S_{1}(1)\oplus S_{2}(1)\oplus S_{3}(1)\oplus S_{3}(1)\oplus S_{4}(1)\oplus S_{5}(1) $\\
&&$\tau^{9}I(3)$&$ 2 $&$        S_{1}'(1)\oplus S_{1}'(1)\oplus S_{2}'(1)\oplus S_{2}'(1)\oplus S_{3}'(1)\oplus S_{3}'(1) $\\
&&$\tau^{14}I(3) $&$ 3 $&$        S_{1}''(1)\oplus S_{1}''(1)\oplus S_{1}''(1)\oplus S_{2}''(1)\oplus S_{2}''(1)\oplus S_{2}''(1) $\\\hline
&$        P(4)        $&$ \tau^{5}I(4)$&$ 1 $&$  S_{1}(1)\oplus S_{2}(1)\oplus S_{3}(1)\oplus S_{4}(1)\oplus S_{5}(1)  $\\
&&$ \tau^{9}I(4) $&$ 2 $&$  S_{1}'(1)\oplus S_{1}'(1)\oplus S_{2}'(1)\oplus S_{3}'(1)\oplus S_{3}'(1)$\\
&&$ \tau^{14}I(4) $&$ 3 $&$  S_{1}''(1)\oplus S_{1}''(1)\oplus S_{1}''(1)\oplus S_{2}''(1)\oplus S_{2}''(1)  $\\\hline
&$        P(5)        $
&$ \tau^{5}I(5)  $&$ 1 $&$  S_{1}(1)\oplus S_{3}(1)\oplus S_{4}(1)\oplus S_{5}(1) $\\
&&$ \tau^{8}I(2) $&$ 1 $&$  S_{1}(1)\oplus S_{2}(2)\oplus S_{4}(1)\oplus S_{5}(1) $\\
&&$ \tau^{8}I(5)  $&$ 1 $&$  S_{1}(1)\oplus S_{2}(2)\oplus S_{3}(2)\oplus S_{5}(1) $\\
&&$ \tau^{11}I(2)  $&$ 1 $&$  S_{1}(1)\oplus S_{2}(2)\oplus S_{3}(2)\oplus S_{4}(2)$\\
&&$ \tau^{9}I(5)  $&$ 2 $&$  S_{1}'(1)\oplus S_{2}'(1)\oplus S_{3}'(1)\oplus S_{3}'(1) $\\
&&$ \tau^{14}I(5)  $&$ 3 $&$  S_{1}''(1)\oplus S_{1}''(1)\oplus S_{2}''(1)\oplus S_{2}''(1) $\\

\hline
&$        P(6)        $
&$ \tau^{5}I(6)  $&$ 1 $&$  S_{2}(1)\oplus S_{3}(1)\oplus S_{4}(1) $\\
&&$ \tau^{9}I(9) $&$ 1 $&$  S_{1}(2)\oplus S_{3}(1)\oplus S_{4}(1) $\\
&&$ \tau^{9}I(6)  $&$ 1 $&$  S_{3}(1)\oplus S_{4}(1)\oplus S_{5}(3) $\\
&&$ \tau^{11}I(9)  $&$ 1 $&$  S_{1}(2)\oplus S_{2}(2)\oplus S_{4}(1) $\\
&&$ \tau^{13}I(9) $&$ 1 $&$  S_{2}(2)\oplus S_{4}(1)\oplus S_{5}(3) $\\
&&$ \tau^{13}I(6)  $&$ 1 $&$  S_{1}(3)\oplus S_{4}(1)\oplus S_{5}(3)$\\
&&$ \tau^{9}I(6)  $&$ 2 $&$  S_{1}'(1)\oplus S_{2}'(1)\oplus S_{3}'(1) $\\
&&$ \tau^{14}I(6)  $&$ 3 $&$  S_{1}''(1)\oplus S_{1}''(1)\oplus S_{2}''(1)$\\

\hline
&$        P(7)        $
&$ \tau^{5}I(7)  $&$ 1 $&$  S_{2}(1)\oplus S_{3}(1) $\\
&&$ \tau^{11}I(1)  $&$ 1 $&$  S_{1}(2)\oplus S_{3}(1)$\\
&&$ \tau^{14}I(1)  $&$ 1 $&$  S_{3}(1)\oplus S_{5}(3)$\\
&&$ \tau^{14}I(7)  $&$ 1 $&$  S_{3}(1)\oplus S_{4}(4)$\\
&&$ \tau^{9}I(7)  $&$ 2 $&$  S_{1}'(1)\oplus S_{2}'(1)$\\
&&$ \tau^{14}I(7)  $&$ 2 $&$  S_{2}'(1)\oplus S_{3}'(2)$\\
&&$ \tau^{14}I(7)  $&$ 3 $&$  S_{1}''(1)\oplus S_{2}''(1)$\\

\hline
&$        P(9)        $
&$ \tau^{5}I(9)  $&$ 1 $&$  S_{2}(1)\oplus S_{4}(1)\oplus S_{5}(1)$\\
&&$ \tau^{7}I(9)  $&$ 1 $&$  S_{2}(1)\oplus S_{3}(2)\oplus S_{5}(1) $\\
&&$ \tau^{11}I(6)  $&$ 1 $&$  S_{2}(1)\oplus S_{3}(2)\oplus S_{4}(2)$ \\
&&$ \tau^{9}I(6) $&$ 1 $&$  S_{1}(2)\oplus S_{4}(1)\oplus S_{5}(1) $\\
&&$ \tau^{9}I(9) $&$ 1 $&$  S_{1}(2)\oplus S_{3}(2)\oplus S_{5}(1) $\\
&&$ \tau^{13}I(6) $&$ 1 $&$  S_{1}(2)\oplus S_{2}(3)\oplus S_{5}(1)$ \\
&&$ \tau^{9}I(9) $&$ 2 $&$  S_{1}'(1)\oplus S_{2}'(1)\oplus S_{3}'(1)$\\
&&$ \tau^{14}I(9)  $&$ 3 $&$  S_{1}''(1)\oplus S_{1}''(1)\oplus S_{2}''(1)$\\

\end{longtable}